\documentclass[a4paper, 12pt, onecolumn]{article}
\usepackage{geometry}

\usepackage{amsmath}
\usepackage{amsthm}
\usepackage{amsfonts}
\usepackage{bbm}
\usepackage{CJK}
\usepackage{fancyhdr}
\usepackage{graphicx}
\usepackage{geometry}
\usepackage{psfrag}
\usepackage{amsfonts,amsmath,amsthm, amssymb}
\usepackage{latexsym, euscript, epic, eepic}
\usepackage{time}
\usepackage{txfonts}
\usepackage{colortbl}
\usepackage{stmaryrd}
\usepackage{mathrsfs}
\usepackage{txfonts}
\usepackage{amsfonts}
\usepackage{color}
\usepackage{lineno}
\usepackage[square, comma, sort&compress, numbers]{natbib}
\usepackage{indentfirst,latexsym,bm}

 \numberwithin{equation}{section}

\begin{document}
\newtheorem{theorem}{Theorem}[section]
\newtheorem{proposition}[theorem]{Proposition}
\newtheorem{remark}[theorem]{Remark}
\newtheorem{corollary}[theorem]{Corollary}
\newtheorem{definition}{Definition}[section]
\newtheorem{lemma}[theorem]{Lemma}
\newtheorem{conjecture}[theorem]{Conjecture}
\newtheorem{question}[theorem]{Question}
\newtheorem{problem}[theorem]{Problem}
\newtheorem{open problem}[theorem]{Open Problem}
\newtheorem{condition}[theorem]{Condition}
\newcommand{\IN}{\mathbb{N}}
\newcommand{\IR}{\mathbb{R}}
\newcommand{\IC}{\mathbb{C}}
\newcommand{\arth}{{\rm{arth}}}
\newcommand{\arsh}{\rm{arsh}}
\newcommand{\K}{\mathscr{K}}
\newcommand{\E}{\mathscr{E}}
\newcommand{\F}{\mathscr{F}}
\newcommand{\M}{\mathscr{M}}

\title{\large
Tur\'an-Type Inequalities for Gaussian\\ Hypergeometric Functions, and Baricz's Conjecture
}
\author{\normalsize Song-Liang Qiu\footnote{{\it Corresponding author.}}\,, Xiao-Yan Ma and Xue-Yan Xiang
}
\date{}
\maketitle
\fontsize{12}{22}\selectfont\small

\paragraph{Abstract:}
In 2007, \'A. Baricz put forward a conjecture concerning Tur\'an-type inequalities for Gaussian hypergeometric functions (see Conjecture \ref{ConjA} in Section \ref{Sec1}). In this paper, the authors disprove this conjecture with several methods, and
present Tur\'an-type double inequalities for Gaussian hypergeometric functions, and sharp bounds for complete and generalized elliptic integrals of the first kind.
\\[10pt]
\\[10pt]

\section{\normalsize Introduction}\label{Sec1}

Throughout this paper, we let $r^{\,\prime}=\sqrt{1-r^2}$ for each $r\in[0,1]$, $\IN$ denote the set of positive integers, $\IN_0=\IN\cup\{0\}$, and let $\Gamma(x)=\int_0^\infty t^{x-1}e^{-t}dt$ , $\psi(x)=\Gamma'(x)/\Gamma(x)$ and $B(x,y)=\Gamma(x)\Gamma(y)/\Gamma(x+y)$ be the classical gamma, psi and beta functions, respectively. For $a\in(0,1)$, we always let 
\begin{align}\label{beta(a)}
\alpha=&\alpha(a)=\frac{3-2a}{3(1-a)}, ~\eta=\eta(a)=\frac{a}{3-2a} \mbox{~ and ~}
\beta=\beta(a)=\frac{\psi(3/2-a)-\psi(1-a)}{B(a,1-a)},
\end{align}
and for $a, b_1, b_2\in(0,\infty)$, $c_1=a+b_1$ and $c_2=a+b_2$ with $b_1<b_2$, let
\begin{align*}
\overline{\alpha}=\overline{\alpha}(a,b_1,b_2)=\frac{b_2c_1}{b_1c_2}, ~\overline{\eta}=\overline{\eta}(a,b_1,b_2)=a\frac{b_2-b_1}{b_2c_1} \mbox{~ and ~}
\overline{\beta}=\overline{\beta}(a,b_1,b_2)=\frac{\psi(b_2)-\psi(b_1)}{B(a,b_1)}.
\end{align*}
Obviously, $\overline{\alpha}(a,1-a,3/2-a)=\alpha(a)$, $\overline{\eta}(a,1-a,3/2-a)=\eta(a)$, $\overline{\beta}(a,1-a,3/2-a)=\beta(a)$,
$$1<\alpha=1+\frac{a}{3(1-a)}, ~0<\eta<a, ~1<\overline{\alpha}=1+\frac{a(b_2-b_1)}{b_1c_2} \mbox{~ and ~} 0<\overline{\eta}<\frac{a}{c_1}.$$

For $a, b, c, x\in\IC$ with $c\neq 0, -1, -2, \cdots$ and $|x|<1$, the Gaussian hypergeometric function is defined by
\begin{align}\label{Def-2F1}
F(a,b;c;x)={}_2F_1(a,b;c;x)=\sum_{n=0}^{\infty}\frac{(a)_n(b)_n}{(c)_nn!}x^n
\end{align}
which is said to be zero-balanced if $c=a+b$, where $(a)_n$ is the Appell (or Pochhammer) symbol defined by $(a)_0=1$ for $a\neq0$, and $(a)_n=a(a+1)(a+2)\cdots(a+n-1)=\Gamma(n+a)/\Gamma(a)$ for $n\in\IN$. As particular cases of $F(a,b;c;x)$, the generalized elliptic integrals of the first and second kinds are defined as
\begin{align}\label{Ka}
\K_a=\K_a(r)=\frac\pi2F\left(a,1-a;1;r^2\right), \,\K_a'=\K_a'(r)=\K_a(r^{\,\prime})
\end{align}
and
\begin{align}\label{Ea}
\E_a=\E_a(r)=\frac\pi2F\left(a-1,1-a;1;r^2\right), \,\E_a'=\E_a'(r)=\E_a(r^{\,\prime})
\end{align}
for $a,r\in(0,1)$, respectively, and the inverse hyperbolic tangent $\arth$ can be represented as
\begin{align}\label{arth}
\arth\,r=rF\left(1/2,1;3/2;r^2\right)
\end{align}
for $r\in[0,1)$ (see \cite{AQVV} and \cite[15.1.4]{AS}). In particular, $\K=\K_{1/2}$ and $\K'=\K_{1/2}'$ \big{(}$\E=\E_{1/2}$ and $\E'=\E_{1/2}'$\big{)} are the well-known complete elliptic integrals of the first (second, respectively) kind.

It is well known that Gaussian hypergeometric functions have many important applications in several fields of mathematics as well as in physics and engineering, and numerous results have been obtained for them. Because of this, revealing their properties has been a main topic in the field of special functions. For the basic properties of Gaussian hypergeometric functions, the readers are referred to \cite{AAR,AQ,AS,Ask2,AVV,AVVb,Be1,Be2,BB3,C,OLB,QMC1,QMC2,QVu}.

One of the topics in the study of hypergeometric functions is to reveal their dependence on parameters (cf. \cite{AQVV,Ba,QMB}, for example).  In \cite{Ba}, \'A. Baricz presented several significant and interesting properties, including Tur\'an-type inequalities, for the zero-balanced hypergeometric functions and some of their particular cases such as generalized elliptic integrals, and put forward the following conjecture.
\begin{conjecture}\label{ConjA}
{\rm (\cite[Conjecture A]{Ba})} For each $a>0$ and $r\in(0,1)$, and for $b\in(0,\infty)$, we have
\begin{align}\label{IneqA}
\frac{\partial}{\partial b}\left[\frac{F(a+1,b+1;a+b+1;r)}{F(a,b;a+b;r)}\right]<0.
\end{align}
In particular, for all $a,r\in(0,1)$ we have the following Tur\'an type inequality
\begin{align}\label{IneqA1}
\frac{F(a+1,2-a;2;r)}{F(a,1-a;1;r)}>\frac{F(a+1,5/2-a;5/2;r)}{F(a,3/2-a;3/2;r)}.
\end{align}
\end{conjecture}
For $a,b\in(0,\infty)$ and $r\in(0,1)$, let $f(a,b,r)=F(a,b;a+b+1;r)/F(a,b;a+b;r)$. Then by the formula
\begin{align}\label{15.3.3}
F(a,b;c;r)=(1-r)^{c-a-b}F(c-a,c-b;c;r)
\end{align}
(cf. \cite[15.3.3]{AS}),  \eqref{IneqA} and \eqref{IneqA1} are equivalent to the following little simpler inequalities
\begin{align}\label{IneqA'}
\frac{\partial f(a,b,r)}{\partial b}=&\frac{\partial}{\partial b}\left[\frac{F(a,b;a+b+1;r)}{F(a,b;a+b;r)}\right]<0
\end{align}
and
\begin{align}
f(a,1-a,r)>f(a,3/2-a,r),\label{IneqA1'}
\end{align}
respectively. Hence in the sequel, we can replace \eqref{IneqA} and \eqref{IneqA1} by \eqref{IneqA'} and \eqref{IneqA1'}, respectively.

For $a,r\in(0,1)$, let $f$ be as in \eqref{IneqA'}, and define $\Lambda$ and $\Lambda_1$ by
\begin{align}
\Lambda(a,r)=\frac{\Lambda_1(a,r)}{F(a,1-a;1;r)}=\frac{f(a,1-a,r)}{f(a,3/2-a,r)}=\frac{F(a,1-a;2;r)F(a,3/2-a;3/2;r)}{F(a,1-a;1;r)F(a,3/2-a;5/2;r)},\label{Lam}
\end{align}
and for $a, b_1, b_2\in(0,\infty)$, $c_1=a+b_1$, $c_2=a+b_2$ and $r\in(0,1)$ with $b_1<b_2$, define $\Lambda_2$ and $\Lambda_3$ by
\begin{align}\label{Lam2}
\Lambda_2(a,b_1,b_2,r)=\frac{\Lambda_3(a,b_1,b_2,r)}{F(a,b_1;c_1;r)}=\frac{f(a,b_1,r)}{f(a,b_2,r)}=\frac{F(a,b_1;c_1+1;r)F(a,b_2;c_2;r)}{F(a,b_1;c_1;r)F(a,b_2;c_2+1;r)}.
\end{align}
The inequalities satisfied by $\Lambda$ and $\Lambda_2$ are Tur\'an-type inequalities. It is easy to see that
\begin{align*}
\Lambda(a,r)=&\Lambda_2(a,1-a,3/2-a,r), ~\Lambda_1(a,r)=\Lambda_3(a,1-a,3/2-a,r),\\
\Lambda(1/2,r)=&\frac{4\arth\,t}{3\K(t)}\cdot\frac{\E(t)-t^{\,\prime2}\K(t)}{t-t^{\,\prime2}\arth\,t} \mbox{~ and ~}
\Lambda_1(1/2,r)=\frac{8\arth\,t}{3\pi}\cdot\frac{\E(t)-t^{\,\prime2}\K(t)}{t-t^{\,\prime2}\arth\,t}
\end{align*}
(cf. \eqref{Form1} in Section 2), where $t=\sqrt{r}$, and by applying the well-known Ramanujan's asymptotic formula for zero-balanced hypergeometric functions (cf. \cite[15.3.10]{AS}), we can easily obtain the limiting values
\begin{align}\label{limLam}
\Lambda_2(a,b_1,b_2,0)=\Lambda_3(a,b_1,b_2,0)=1 \mbox{~ and ~} \Lambda_2(a,b_1,b_2,1^-)=\overline{\alpha}.
\end{align}
In particular, $\Lambda(a,0)=\Lambda_1(a,0)=1$ and $\Lambda(a,1^-)=\alpha$. It is clear that
\begin{align}\label{Lam23}
\eqref{IneqA'} \mbox{~ holds~} \Longleftrightarrow\Lambda_2(a,b_1,b_2,r)>1\Longleftrightarrow\Lambda_3(a,b_1,b_2,r)>F(a,b_1;c_1;r),
\end{align}
and in particular,
\begin{align}\label{Lam01}
\eqref{IneqA1'} \mbox{~ holds~} \Longleftrightarrow\Lambda(a,r)>1\Longleftrightarrow\Lambda_1(a,r)>F(a,1-a;1;r).
\end{align}

So far, there is no proof or disproof of Conjecture \ref{ConjA} in the literature yet. By \eqref{limLam}--\eqref{Lam01}, it is natural to think of that a method for us to study Conjecture \ref{ConjA} is to reveal the properties of the following functions:
\begin{align}
f_1(a,r)\equiv&\Lambda_1(a,r)-F(a,1-a;1;r),\label{f1}\\
f_2(a,r)\equiv&\alpha F(a,1-a;1;r)-\Lambda_1(a,r),\label{f2}\\
f_3(a,b_1,b_2,r)\equiv&\Lambda_3(a,b_1,b_2,r)-F(a,b_1;a+b_1;r),\label{f3}\\
f_4(a,b_1,b_2,r)\equiv&\overline{\alpha}F(a,b_1;a+b_1;r)-\Lambda_3(a,b_1,b_2,r).\label{f4}
\end{align}
Obviously, $f_3(a,1-a,3/2-a,r)=f_1(a,r)$, $f_4(a,1-a,3/2-a,r)=f_2(a,r)$, $f_1(a,0)=f_3(a,b_1,b_2,0)=0$, and \eqref{Lam23} (\eqref{Lam01}) is true if and only if $f_3$ ($f_1$, respectively) is a positive function. On the other hand, the following question is naturally raised.
\begin{question}\label{Q}
Is it true that $f_2(a,r)>0$ ($f_4(a,b_1,b_2,r)>0$) for all $a,r\in(0,1)$ ($a, b_1, b_2\in(0,\infty)$ and $r\in(0,1)$ with $b_1<b_2$, respectively)?
\end{question}
If the answer to Question \ref{Q} is affirmative, and if Conjecture \ref{ConjA} is true, then the Tur\'an-type double inequality
$1<\Lambda(a,r)<\alpha$ $\left(1<\Lambda_2(a,b_1,b_2,r)<\overline{\alpha}\right)$ holds for $a,r\in(0,1)$ ($a, b_1, b_2\in(0,\infty)$ and $r\in(0,1)$ with $b_1<b_2$, respectively).

The main purpose of this paper is to disprove Conjecture \ref{ConjA}, give an affirmative answer to Question \ref{Q}, and present Tur\'an-type double inequalities for $\Lambda(a,r)$ and $\Lambda_2(a,b_1,b_2,r)$, mainly by studying the properties of $f_j$ ($1\leq j\leq4$). These results are stated below as Theorems \ref{th1.1} and \ref{th1.2}, whose proofs will be given in Section \ref{Sec3} and Section \ref{Sec5}, respectively. In order to prove these two theorems, we shall present several preliminary results in Section \ref{Sec2}. In Section \ref{Sec4} (Section \ref{Sec6}), we shall give a detailed counterexample for Conjecture \ref{ConjA} (remarks giving some other methods to disprove Conjecture \ref{ConjA}, three open problems and a conjecture, respectively).

Throughout this paper, we always use the notations introduced above, without repeating their definitions.

\begin{theorem}\label{th1.1}
(1) Conjecture \ref{ConjA} is not true if $b>\sqrt{a(a+1)}$. In particular, \eqref{IneqA1} does not hold if $a\in(0,1/3)$.

(2) If
$a\in(0,3/7)$ ($a\in(3/7,1)$), then there exists a number $r_1\in(0,1)$ ($r_2\in(0,1)$) such that $f_1$ is strictly decreasing (increasing) in $r\in(0,r_1]$ ($r\in(0,r_2]$), so that $\Lambda(a,r)<1 ~~(\Lambda(a,r)>1)$
for $r\in(0,r_1]$ ($r\in(0,r_2]$, respectively). Moreover, for each $a\in(0,1)$, there exists a number $r_3\in(0,1)$ such that $f_1$ is positive and strictly increasing in $r\in[r_3,1)$ with $f_1(a,1^-)=\infty$ so that $\Lambda(a,r)>1$
for $r\in[r_3,1)$. In particular, neither the first (or equivalently, the second) inequality in \eqref{Lam01} nor its reversed one is always valid for all $a,r\in(0,1)$.

(3) For $n\in\IN_0$ and $a, r\in(0,1)$, suppose that $f_2(a,r)=\sum_{k=0}^\infty a_kr^k$, $P_n(a,r)=\sum_{k=0}^{n+1}a_kr^{2k}$ and $\delta_n=\alpha\beta-\sum_{k=0}^{n+1}a_k$, where $a_k=a_k(a)$ depends on $a$, and let $a_k^{(1)}=a_k^{(1)}(a,1-a,3/2-a)$ be as in Lemma \ref{lm2.3}. Then $f_2(a,0)=a_0=\alpha-1=a/[3(1-a)]$, $f_2(a,1^-)=\alpha\beta$, and
\begin{align}\label{an}
a_k=a_k(a)=\frac{1}{k}\sum_{i=0}^{k-1}\frac{(i+1)(a)_{k-i-1}(3/2-a)_{k-i-1}}{(3/2)_{k-i-1}(k-i-1)!}\big{|}a_{i+1}^{(1)}\big{|}>0
\end{align}
for $k\in\IN$, with $a_1=a(a+1)/10$, $a_2=a(a+1)(-41a^2+24a+100)/2100$, and $a_k(0^+)=0$ for $k\in\IN_0$.
In particular, for each $a\in(0,1)$, $f_2$ is absolutely monotone in $r$ from $(0,1)$ onto $(a/[3(1-a)],\alpha\beta)$, $\delta_n$ is strictly decreasing to $0$ in $n\in\IN_0$, and for all $a, r\in(0,1)$,
\begin{align}\label{Expr1Ka}
\alpha\K_a(r)=\frac{\E_a(r)-r^{\,\prime 2}\K_a(r)}{ar^2}\cdot\frac{F\bigl(a,3/2-a;3/2;r^2\bigr)}{F\bigl(a,3/2-a;5/2;r^2\bigr)}+\frac\pi2\sum_{k=0}^\infty a_kr^{2k}
\end{align}
and
\begin{align}\label{Expr1K}
\K(r)=&\frac{\E(r)-r^{\,\prime 2}\K(r)}{1-\left(r^{\,\prime 2}\arth\,r\right)/r}\frac{\arth\,r}{r}+\frac{3\pi}{8}\sum_{k=0}^\infty a_k(1/2)r^{2k}\nonumber\\
=&\frac{3\pi}{8}\left[\left(\sum_{k=0}^\infty a_k^{(2)}r^{2k}\right)\frac{\arth\,r}{r}+\sum_{k=0}^\infty a_k(1/2)r^{2k}\right],
\end{align}
where $a_k^{(2)}$ is defined by \eqref{Defan} in Corollary \ref{Col2.4}. Moreover,
\begin{align}
0<&\frac{2\alpha}{\pi}\K_a(r)-2\frac{\E_a(r)-r^{\,\prime 2}\K_a(r)}{ \pi ar^2}\cdot\frac{F\bigl(a,3/2-a;3/2;r^2\bigr)}{F\bigl(a,3/2-a;5/2;r^2\bigr)}-P_n(a,r)<\delta_n r^{2(n+2)},\label{Ineqf2}\\
a_{n+2}r^{n+2}<&\alpha F(a,1-a;1;r)-P_n\left(a,\sqrt{r}\right)-\Lambda_1(a,r)<\delta_n r^{n+2},\label{IneqLam1}\\
-\delta_nr^{n+1}<&\frac{\Lambda_1(a,r)-F(a,1-a;1;r)}{r}-\frac{a(7a-3)}{30}-Q_n(a,r)<-a_{n+2}r^{n+1}\label{IneqLam1'}
\end{align}
and
\begin{align}\label{IneqA1''}
\alpha(1-\beta)<&\alpha\left[1-\frac{\beta}{F(a,1-a;1;r)}\right]<\alpha-\frac{\alpha\beta r^{n+2}}{F(a,1-a;1;r)}<\alpha-\frac{P_n\bigl(a,\sqrt{r}\bigr)+\delta_n r^{n+2}}{F(a,1-a;1;r)}<\Lambda(a,r)\nonumber\\
<&\alpha-\frac{P_{n+1}\bigl(a,\sqrt{r}\bigr)}{F(a,1-a;1;r)}<
\alpha-\frac{r^{n+2}\sum_{k=0}^{n+2}a_k+a_0\bigl(1-r^{n+2}\bigr)}{F(a,1-a;1;r)}<\alpha-\frac{\eta}{F(a,1-a;1;r)}<\alpha
\end{align}
for each $n\in\IN_0$ and for all $a, r\in(0,1)$, where
\begin{align*}
Q_n(a,r)=&\frac{a}{3(1-a)}\sum_{k=1}^\infty\frac{(a)_{k+1}(1-a)_{k+1}}{[(k+1)!]^2}r^k
-\sum_{k=1}^na_{k+1}r^k\\
=&a\frac{F(a,1-a;1;r)-1-a(1-a)r}{3(1-a)r}-\sum_{k=1}^na_{k+1}r^k.
\end{align*}
\eqref{Ineqf2} is sharp as $r\to0$ (or $r\to1$ or $a\to0$), each inequality in \eqref{IneqA1''} is sharp as $a\to0$, and as constant depending only on $a$, each $\alpha$ in \eqref{Ineqf2}--\eqref{IneqLam1} and  \eqref{IneqA1''} is best possible.
\end{theorem}

\begin{theorem}\label{th1.2}
(1) For $a, b_1, b_2\in(0,\infty)$ and $r\in(0,1)$ with $b_1<b_2$, each of the two inequalities in \eqref{Lam23} is not always valid, provided that $b_1b_2>a(a+1)$. In particular, Conjecture \ref{ConjA} is not true.

(2) For $a\in(0,\infty)$ and $b_1, b_2\in(0,\infty)$ with $b_1<b_2$, if $b_1b_2>a(a+1)$ ($b_1b_2<a(a+1)$), then there exists a number $r_4\in(0,1)$ ($r_5\in(0,1)$) such that $\Lambda_2(a,b_1,b_2,r)<1$ ($\Lambda_2(a,b_1,b_2,r)>1$)
for $r\in(0,r_4)$ ($r\in(0,r_5)$, respectively).
Moreover, for $a,b_1,b_2\in(0,\infty)$ with $b_1<b_2$, there exists a number $r_6\in(0,1)$ such that $\Lambda_2(a,b_1,b_2,r)>1$ for $r\in(r_6,1)$.
In particular, neither the first (or equivalently, second) inequality in \eqref{Lam23} nor its reverse inequality always holds for all $a, b_1, b_2\in(0,\infty)$ and $r\in(0,1)$ with $b_1<b_2$.

(3) The function $f_4$ is strictly increasing in $r$ from $(0,1)$ onto $(a(b_2-b_1)/(b_1c_2),\overline{\alpha}\overline{\beta})$. In particular, for $a, b_1, b_2\in(0,\infty)$ and $r\in(0,1)$ with $b_1<b_2$,
\begin{align}
\overline{\alpha}\left[F(a,b_1;c_1;r)-\overline{\beta}\right]&\leq\Lambda_3(a,b_1,b_2,r)
\leq\overline{\alpha}\left[F(a,b_1;c_1;r)-\overline{\eta}\right],\label{IneqLam2}\\
\overline{\alpha}\left[1-\frac{\overline{\beta}}{F(a,b_1;c_1;r)}\right]
&\leq\Lambda_2(a,b_1,b_2,r)\leq\overline{\alpha}\left[1-\frac{\overline{\eta}}{F(a,b_1;c_1;r)}\right]\label{IneqLam3}
\end{align}
and
\begin{align}
\overline{\alpha}(1-\overline{\beta})&\leq\Lambda_2(a,b_1,b_2,r)\leq\overline{\alpha},\label{Ineq2g1}
\end{align}
with equality in each instance if and only if $b_1=b_1$. Moreover, as constants independent of $r$, the coefficients in the lower and upper bounds in \eqref{IneqLam2}--\eqref{IneqLam3} and the upper bound in \eqref{Ineq2g1} are all best possible, \eqref{IneqLam3} is sharp as $r\to1^-$, and \eqref{IneqLam2}--\eqref{Ineq2g1} are all sharp as $a\to0^+$.

(4) For $n\in\IN_0$, $a, b_1, b_2\in(0,\infty)$, $c_1=a+b_1$, $c_2=a+b_2$ with $b_1<b_2$, and for $r\in(0,1)$, suppose that $f_4(a,b_1,b_2,r)=\sum_{k=0}^\infty \tilde{b}_kr^k$, $\eta_n=\overline{\alpha}\overline{\beta}-\sum_{k=0}^{n+1}\tilde{b}_k$, where $\tilde{b}_k=\tilde{b}_k(a,b_1,b_2)$ depends on $a, b_1$ and $b_2$, and let
\begin{align}\label{Barb}
\bar{b}=\frac{\sqrt{a^2+4a-1}-\sqrt{a^2-1}}{2}\sqrt{\frac{a+1}{a-1}} ~~(>1)
\end{align}
for $a\in(1,\infty)$. Then $f_4$ is absolutely monotone in $r$ from $(0,1)$ onto $(a(b_2-b_1)/(b_1c_2),\overline{\alpha}\overline{\beta})$, provided that $a\in(0,1]$, or $a\in(1,\infty)$ and $b_2\in(0,\bar{b}]$. In particular, for each $n\in\IN_0$ and all $r\in(0,1)$,
\begin{align}\label{IneqLam23}
\sum_{k=0}^{n+2}\tilde{b}_kr^k<&\overline{\alpha}F\left(a,b_1;c_1;r\right)-
\Lambda_3(a,b_1,b_2,r)<\sum_{k=0}^{n+1}\tilde{b}_kr^k+\eta_nr^{n+2}
\end{align}
and
\begin{align}\label{IneqLam23'}
\overline{\alpha}(1-\overline{\beta})<\overline{\alpha}-\frac{\sum_{k=0}^{n+1}\tilde{b}_kr^k+\eta_nr^{n+2}}{F\left(a,b_1;c_1;r\right)}<&\Lambda_2(a,b_1,b_2,r)
<\overline{\alpha}-\frac{\sum_{k=0}^{n+2}\tilde{b}_kr^k}{F\left(a,b_1;c_1;r\right)}<\overline{\alpha}.
\end{align}
\end{theorem}

In the sequel, for $n\in\IN_0$ and $a, b\in(0,\infty)$, we always let
\begin{align*}
\varphi_n(b)=&\frac{(b)_n}{(a+b+1)_n}, ~\rho_n(b)=\frac{(b)_n}{(a+b)_n},\\
\tau_n(b)=&\psi(a+b+1)-\psi(b)+\psi(n+b)-\psi(n+a+b+1),\\
\lambda_n(b)=&\psi(a+b)-\psi(b)+\psi(n+b)-\psi(n+a+b),
\end{align*}
$R(a,b)=-2\gamma-\psi(a)-\psi(b)$ which is called the Ramanujan constant, and let $R(a)=R(a,1-a)$ for $a\in(0,1)$.
By the basic properties of $\Gamma$ and $\psi$, one can easily obtain the following relations and values:
\begin{align}
\varphi_n(b)=&\frac{c}{n+c}\rho_n(b), \varphi_{n+1}(b)=\frac{n+b}{n+c+1}\varphi_n(b),\label{vrn+1}\\
\rho_{n+1}(b)=&\frac{n+b}{n+c}\rho_n(b), \tau_n(b)=\lambda_n(b)+\frac{n}{c(n+c)},\label{taulam}\\
\varphi_0(b)-1=&\rho_0(b)-1=\tau_0(b)=\lambda_0(b)=0, \varphi_1(b)=\frac{b}{c+1}, \rho_1(b)=\frac{b}{c},\label{v01ro1}
\end{align}
and for each $n\in\IN$,
\begin{align}
\varphi_n(0^+)=&\rho_n(0^+)=\varphi_n(\infty)-1=\rho_n(\infty)-1=\varphi_\infty(b)=\rho_\infty(b)=0.\label{vrn0infty}
\end{align}

\section{\normalsize Preliminaries}\label{Sec2}

In this section, we prove several preliminary results needed in the proofs of Theorems \ref{th1.1} and \ref{th1.2}. First, we prove the following Lemma \ref{lm2.1} by which one can easily verify that in what follows, all the processes of calculating limits and derivatives for series of functions term by term are feasible.

\begin{lemma}\label{lm2.1}
(1) For $a, b\in(0,\infty)$ ($n\in\IN$) and $c=a+b$, $\tau_n$ and $\lambda_n$ are both strictly increasing (decreasing) in $n\in\IN$ (in $b\in(0,\infty)$, respectively). Moreover, for $n\in\IN$ and $a,b\in(0,\infty)$,
$\tau_\infty(b)=\psi(c+1)-\psi(b)$, $\lambda_\infty(b)=\psi(c)-\psi(b)$, $\tau_n(0^+)=\lambda_n(0^+)=\infty$, $\tau_n(\infty)=\lambda_n(\infty)=\lim_{b\to\infty}[\varphi_n(b)\tau_n(b)]
=\lim_{b\to\infty}[\rho_n(b)\lambda_n(b)]=0$,
\begin{align*}
\lim_{b\to0^+}[\varphi_n(b)\tau_n(b)]=&\frac{(n-1)!}{(a+1)_n} \mbox{~ and ~}
\lim_{b\to0^+}[\rho_n(b)\lambda_n(b)]=\frac{(n-1)!}{(a)_n}.
\end{align*}

(2) The sequences $\{\varphi_n(b)\}_{n=0}^\infty$ and $\{\rho_n(b)\}_{n=0}^\infty$ are both strictly deceasing, and for each $n\in\IN$, $\varphi_n(b)$ and $\rho_n(b)$ are both strictly increasing and log-concave in $b$ from $(0,\infty)$ onto $(0,1)$.

(3) For $b\in\bigl(0,\sqrt{(a+1)(a+2)}\bigl]$, the sequence $\{\varphi_n(b)\tau_n(b)\}_{n=1}^\infty$ is strictly decreasing in $n\in\IN$. For  $b\in\bigl(\sqrt{(a+1)(a+2)},\infty\bigl)$, there exists an $n_1\in\IN$ such that $\{\varphi_n(b)\tau_n(b)\}_{n=1}^\infty$ is strictly increasing in $n\in\{1, 2, ..., n_1\}$ and then decreasing in $n\geq n_1$.

(4) If $b\in\bigl(0,\sqrt{a(a+1)}\bigl]$, then the sequence $\{\rho_n(b)\lambda_n(b)\}_{n=1}^\infty$ is strictly decreasing in $n\in\IN$. For $b\in\bigl(\sqrt{a(a+1)},\infty\bigl)$, there exists an $n_2\in\IN$ such that $\{\rho_n(b)\lambda_n(b)\}_{n=1}^\infty$ is strictly increasing in $n\in\{1, 2, ..., n_2\}$ and decreasing in $n\geq n_2$.
\end{lemma}

\begin{proof}
(1) It follows from \cite[6.3.5 \& 6.3.16]{AS} that for $n\in\IN$,
\begin{align}\label{taun}
\tau_n(b)=&\sum_{k=0}^\infty\frac{a+1}{(k+b)(k+c+1)}-\sum_{k=0}^\infty\frac{a+1}{(k+n+b)(k+n+c+1)}\nonumber\\
=&\sum_{k=0}^{n-1}\frac{a+1}{(k+b)(k+c+1)}.
\end{align}
Similarly,
\begin{align}\label{lamn'}
\lambda_n(b)=\sum_{k=0}^{n-1}\frac{a}{(k+b)(k+c)}.
\end{align}
The first assertion in part (1) now follows from \eqref{taun} and \eqref{lamn'}.

 Applying \cite[6.1.47, 6.3.5 \& 6.3.18]{AS} and \eqref{v01ro1}--\eqref{vrn0infty}, we obtain the following limiting values
\begin{align}
\tau_\infty(b)=&\lim_{n\to\infty}[\psi(c+1)-\psi(b)+\psi(n+b)-\psi(n+c+1)]=\psi(c+1)-\psi(b),\label{tinf}\\
\tau_n(0^+)=&\lim_{b\to0^+}\left[\psi(c+1)-\psi(b+1)+\frac1b+\psi(n+b)-\psi(n+c+1)\right]=\infty,\label{t(0)}\\
\lambda_\infty(b)=&\lim_{n\to\infty}[\psi(c)-\psi(b)+\psi(n+b)-\psi(n+c)]=\psi(c)-\psi(b),\label{linf}\\
\lambda_n(0^+)=&\lim_{b\to0^+}\left[\psi(c)-\psi(b+1)+\frac1b+\psi(n+b)-\psi(n+c)\right]=\infty,\label{l(0)}\\
\tau_n(\infty)=&\lambda_n(\infty)=\lim_{b\to\infty}[\varphi_n(b)\tau_n(b)]
=\lim_{b\to\infty}[\rho_n(b)\lambda_n(b)]=0,\label{tlinf}\\
\lim_{b\to0^+}[\varphi_n(b)\tau_n(b)]=&\lim_{b\to0^+}\frac{b\Gamma(c+1)\Gamma(n+b)}{\Gamma(b+1)\Gamma(n+c+1)}
[\psi(c+1)-\psi(b+1)+1/b\nonumber\\
&+\psi(n+b)-\psi(n+c+1)]=\frac{(n-1)!}{(a+1)_n}\label{vt(0)}
\end{align}
and
\begin{align}
\lim_{b\to0^+}[\rho_n(b)\lambda_n(b)]=&\lim_{b\to0^+}\frac{b\Gamma(c)\Gamma(n+b)[\psi(c)-\psi(b+1)+1/b+\psi(n+b)-\psi(n+c)]}{\Gamma(b+1)\Gamma(n+c)}
=\frac{(n-1)!}{(a)_n}.\label{rl(0)}
\end{align}

(2) The first assertion in part (2) follows from \eqref{vrn+1} and \eqref{taulam}.

By logarithmic differentiation, $\varphi_n'(b)/\varphi_n(b)=\tau_n(b)$ and $\rho_n'(b)/\rho_n(b)=\lambda_n(b)$
for $n\in\IN$, and $\varphi_0'(b)=\rho_0'(b)=0$. Hence the second assertion in part (2) follows from part (1) and \eqref{vrn0infty}.

(3) It follows from \eqref{vrn+1} and \eqref{taun} that for $n\in\IN$ and $c=a+b$,
\begin{align*}
\frac{n+c+1}{(a+1)\varphi_n(b)}\left[\varphi_{n+1}(b)\tau_{n+1}(b)-\varphi_n(b)\tau_n(b)\right]=\frac{1}{n+c+1}-\tau_n(b),
\end{align*}
which is strictly decreasing in $n\in\IN$ by part (1) with
\begin{align*}
\lim_{n\to\infty}\left[\frac{1}{n+c+1}-\tau_n(b)\right]=&\psi(b)-\psi(c+1)<0, \mbox{~ and}\\
\frac{1}{n+c+1}-\tau_n(b)=&\frac{b^2-(a+1)(a+2)}{b(c+1)(c+2)}
\end{align*}
for $n=1$.
Hence part (3) follows.

(4) The proof of part (4) is similar to that of part (3).
\end{proof}

\begin{lemma}\label{lm2.2}
(1) As a function of $a$, $\beta$ defined by \eqref{beta(a)} is strictly increasing from $(0,1)$ onto itself, with $\beta(1/2)=(\log4)/\pi$. Moreover, $\lim_{a\to1^-}[1-\beta(a)]/(1-a)=\log 4$.

(2) For each $r\in(0,1)$, the functions
$h_1(a)\equiv1-\beta(a)/F(a,1-a;1;r)$ and $h_2(a)\equiv1-\eta(a)/F(a,1-a;1;r)$ are both strictly decreasing from $(0,1)$ onto itself. Moreover,  $\lim_{a\to1^-}h_1(a)/(1-a)=\log(4/(1-r))$ and $\lim_{a\to1^-}h_2(a)/(1-a)=3-\log(1-r)$.
\end{lemma}

\begin{proof}
(1) Clearly, $\beta(0^+)=0$, $\beta(1/2)=(\log4)/\pi$ and $\beta(1^-)=1$. Let $h_3(a)=\psi(3/2-a)-\psi(1-a)$ for $a\in(0,1)$. Then it follows from the well-known monotonicity of $\psi^{(n)}$ (cf. \cite[6.4.1 or 6.4.10]{AS}) that
$$h_3^{(n)}(a)=(-1)^n\left[\psi^{(n)}(3/2-a)-\psi^{(n)}(1-a)\right]>0$$
for all $n\in\IN_0$ and $a\in(0,1)$, and hence $h_3$ is absolutely monotone from $(0,1)$ onto $(2-\log4,\infty)$. Since $\beta(a)=h_3(a)[\sin(\pi a)]/\pi$ by the reflection formula (cf. \cite[6.1.17]{AS})
\begin{align}\label{rf}
\Gamma(a)\Gamma(1-a)=\pi/\sin(\pi a),
\end{align}
$\beta$ is strictly increasing from $(0,1/2]$ onto $(0,(\log4)/\pi)$.

Next, assume that $a\in(1/2,1)$. Set $h_4(a)=(1-a)[\psi(2-a)-\psi(3/2-a)]$ and $h_5(a)=\Gamma(a)\Gamma(2-a)$. Then $h_4(1^-)=0$, $h_4(1/2)=1-\log2>0$ by \cite[6.3.4]{AS}, and
\begin{align}\label{h1}
\beta(a)=&\frac{1-a}{\Gamma(a)\Gamma(2-a)}\left[\psi\left(\frac32-a\right)-\psi(2-a)+\frac{1}{1-a}\right]
=\frac{1-h_4(a)}{h_5(a)}
\end{align}
by \cite[6.3.5]{AS}. It is easy to see that $h_5$ is strictly decreasing and convex from $(0,1)$ onto $(1,\infty)$, and
\begin{align}\label{SerOfh5}
h_4(a)=&(1-a)\left[\frac{2}{3-2a}-\frac{1}{2-a}+\sum_{n=1}^\infty\frac1n\left(\frac{2-a}{n+2-a}-\frac{3-2a}{2n+3-2a}\right)\right]\nonumber\\
=&\sum_{n=0}^\infty\frac{1-a}{(n+2-a)(2n+3-2a)}
\end{align}
by \cite[6.3.16]{AS}. It is clear that for each $n\in\IN_0$, the function $a\mapsto (1-a)/[(n+2-a)(2n+3-2a)]$
is strictly decreasing on $[1/2,1)$, and so is $h_4$ by \eqref{SerOfh5}. Hence $1-h_4(a)$ is strictly increasing from $[1/2,1)$ onto $[\log2,1)$, and it follows from \eqref{h1} that $\beta$ is strictly increasing from $(1/2,1)$ onto $((\log4)/\pi,1)$.

The monotonicity of $\beta$ on $(0,1)$ now follows from the above investigation.

By \eqref{rf}, \cite[6.3.5 \& 6.4.6]{AS} and by differentiation,
\begin{align}\label{h1'}
\beta'(a)=&\frac{\sin(\pi a)}{\pi}\left[\psi'(1-a)-\psi'\left(\frac32-a\right)\right]+\left[\psi\left(\frac32-a\right)-\psi(1-a)\right]\cos(\pi a)\nonumber\\
=&\frac{\sin(\pi a)}{\pi}\left[\psi'(2-a)-\psi'\left(\frac32-a\right)\right]+\frac{\sin(\pi a)}{\pi(1-a)^2}\nonumber\\
&+\left[\psi\left(\frac32-a\right)-\psi(2-a)\right]\cos(\pi a)+\frac{\cos(\pi a)}{1-a}.
\end{align}
Applying \eqref {h1'} and \cite[6.3.2 \& 6.3.3]{AS}, we obtain the limiting values $\beta'(0^+)=2(1-\log2)$ and
\begin{align*}
\beta'(1^-)=&\psi(1)-\psi(1/2)+\lim_{a\to1^-}\frac{\sin(\pi a)+\pi(1-a)\cos(\pi a)}{\pi(1-a)^2}
=\psi(1)-\psi(1/2)=\log4.
\end{align*}
Hence by l'H\^opital's rule, $\lim_{a\to1^-}[1-\beta(a)]/(1-a)=\beta'(1^-)=\log4$.

(2) For each $n\in\IN_0$, let $h_{6,\,n}(a)\equiv\Gamma(n+a)\Gamma(n+1-a)$. Then it is easy to see that $h_{6,\,n}$ is strictly decreasing (increasing) on $(0,1/2]$ ( $[1/2,1)$, respectively) with $h_{6,\,n}(0^+)=h_{6,\,n}(1^-)=n!(n-1)!$ for $n\in\IN$, and
\begin{align}\label{SerOfF}
1+\frac{a(1-a)}{\Gamma(a+1)\Gamma(2-a)}\sum_{n=1}^\infty\frac{r^n}{n^2}&<F(a,1-a;1;r) \nonumber\\
=1+\frac{a(1-a)}{\Gamma(a+1)\Gamma(2-a)}\sum_{n=1}^\infty\frac{h_{6,\,n}(a)}{(n!)^2}r^n&<1+\frac{a(1-a)}{\Gamma(a+1)\Gamma(2-a)}\sum_{n=1}^\infty\frac{r^n}{n}
\end{align}
by which we obtain the limiting values
\begin{align}\label{F0}
\lim_{a\to0^+}F(a,1-a;1;r)=\lim_{a\to1^-}F(a,1-a;1;r)=1
\end{align}
(see also \cite[Corollary 7.3(2)]{AQVV}).  Clearly, $h_1(1^-)=0$ by part (1) and \eqref{F0}, and $h_1(a)$ can be written as
$$h_1(a)=1-a\frac{\psi(3/2-a)-\psi(2-a)+1/(1-a)}{\Gamma(a+1)\Gamma(1-a)F(a,1-a;1;r)},$$
by which we obtain the limiting value $h_1(0^+)=1$.

Let $h_7(a)=\Gamma(a)\Gamma(1-a)F(a,1-a;1;r)$, which has the following series expansion
\begin{align*}
h_7(a)=\frac{\pi}{\sin(\pi a)}+\sum_{n=1}^\infty\frac{h_{6,\,n}(a)}{(n!)^2}r^n
\end{align*}
by \eqref{rf} and \eqref{SerOfF}. Obviously, if $a\in(0,1/2]$, then $h_7$ is strictly decreasing on $(0,1/2]$, with $h_7(0^+)=\infty$ and $h_7(1/2)=2\K\bigl(\sqrt{r}\bigr)$. Since $h_1(a)=1-h_3(a)/h_7(a)$, $h_1$ is strictly decreasing on $(0,1/2]$ by the monotonicity properties of $h_3$ and $h_7$.

For $a\in[1/2,1)$, let $h_8(a)=h_5(a)F(a,1-a;1;r)$. Then $h_8(1^-)=1$ by \eqref{F0}, and
\begin{align}\label{h8}
h_1(a)=1-(1-a)\frac{h_3(a)}{h_8(a)}=1-\frac{1-h_4(a)}{h_8(a)}
\end{align}
by \cite[6.3.5]{AS}. By the symmetry of the parameter and \cite[Corollary 7.3(2)]{AQVV}, $F(a,1-a;1;r)$ is strictly decreasing in $a\in[1/2,1)$ for each $r\in(0,1)$, and hence so is $h_8$ by the monotonicity of $h_5$. Since the function $a\mapsto 1-h_4(a)$ is strictly increasing from $[1/2,1)$ onto $[\log2,1)$, $h_1$ is strictly decreasing on $[1/2,1)$ by \eqref{h8}.

From the above discussion, we obtain the monotonicity of $h_1$ on $(0,1)$.

By differentiation, \cite[6.3.5]{AS} and \eqref{SerOfF}, and by Lemma \ref{lm2.1},
\begin{align}\label{partialF}
\frac{\partial F(a,1-a;1;r)}{\partial a}&=\sum_{n=1}^\infty\frac{(a)_n(1-a)_n}{(n!)^2}[\psi(n+a)-\psi(n+1-a)-\psi(a)+\psi(1-a)]r^n\nonumber\\
&=\sum_{n=1}^\infty\frac{a(1-a)[\psi(n+a)-\psi(n+1-a)-\psi(a+1)+\psi(2-a)]+1-2a}{\Gamma(a+1)\Gamma(2-a)(n!)^2}h_{6,\,n}(a)r^n,
\end{align}
which is negative for $a\in[1/2,1)$ by \cite[Theorem 2.5]{Ba}, with
\begin{align}\label{F'01}
&\lim_{a\to0^+}\frac{\partial F(a,1-a;1;r)}{\partial a}=-\lim_{a\to1^-}\frac{\partial F(a,1-a;1;r)}{\partial a}
=\sum_{n=1}^\infty\frac{r^n}{n}=\log\frac{1}{1-r}.
\end{align}
It follows from \eqref{F0} and \eqref{F'01} that
\begin{align*}
h_8'(1^-)=&\lim_{a\to1^-}h_5(a)\left\{F(a,1-a;1;r)[\psi(a)-\psi(2-a)]+\frac{\partial F(a,1-a;1;r)}{\partial a}\right\}\\
=&\lim_{a\to1^-}\frac{\partial F(a,1-a;1;r)}{\partial a}=\log(1-r),
\end{align*}
and hence by \eqref{h8} and l'H\^opital's rule,
\begin{align*}
\lim_{a\to1^-}\frac{h_1(a)}{1-a}&=\lim_{a\to1^-}\frac{h_8(a)+h_4(a)-1}{(1-a) h_8(a)}=\lim_{a\to1^-}\frac{h_8(a)+h_4(a)-1}{1-a}\\
&=\lim_{a\to1^-}\left[\psi(2-a)-\psi\left(\frac32-a\right)+\frac{h_8(a)-1}{1-a}\right]\\
&=\log4+\lim_{a\to1^-}\frac{h_8(a)-1}{1-a}=\log4-h_8'(1^-)=\log\frac{4}{1-r}.
\end{align*}

Clearly, $h_2(0^+)=1$, and $h_2(1^-)=0$ by \eqref{F0}. By l'H\^opital's rule, \eqref{F0} and \eqref{F'01}, we obtain
\begin{align*}
\lim_{a\to1^-}\frac{h_2(a)}{1-a}=&\lim_{a\to1^-}\frac{(3-2a)F(a,1-a;1;r)-a}{(1-a)(3-2a)F(a,1-a;1;r)}\\
=&\lim_{a\to1^-}\left[3F(a,1-a;1;r)+a\frac{F(a,1-a;1;r)-1}{1-a}\right]\\
=&3-\lim_{a\to1^-}\frac{\partial F(a,1-a;1;r)}{\partial a}=3-\log(1-r).
\end{align*}

By \eqref{Ka}, $h_2(a)$ can be written as
\begin{align*}
h_2(a)=&1-\frac{\pi a}{2(3-2a)[\pi F(a,1-a;1;r)/2]}=1-\frac\pi2\left[(3-2a)\cdot\frac{\K_a\bigl(\sqrt{r}\bigr)}{a}\right]^{-1},
\end{align*}
and hence $h_2$ is strictly decreasing on $(0,1/2]$ by \cite[Corollary 7.3(2)]{AQVV}) and \cite[Theorem 1.1(1)]{QMB}.

Suppose that $a\in[1/2,1)$. Then $F(a,1-a;1;r)$ is strictly decreasing in $a\in[1/2,1)$ for each given $r\in(0,1)$ as above-mentioned, and so is the function $a\mapsto(3-2a)F(a,1-a;1;r)$. Hence $h_2$ is strictly decreasing on $[1/2,1)$, too. Consequently, $h_2$ is strictly decreasing on $(0,1)$ by the above investigation.
\end{proof}

\begin{lemma}\label{lm2.3}
(1) For $a, b_1, b_2\in(0,\infty), c_1=a+b_1$ and $c_2=a+b_2$ with $b_1<b_2$, and for $r\in(0,1)$, suppose that
$h_9(r)=F(a,b_1;c_1+1;r)/F(a,b_2;c_2+1;r)=\sum_{n=0}^\infty a_n^{(1)}r^n$,
where $a_n^{(1)}=a_n^{(1)}(a,b_1,b_2)$ depends on $a, b_1$ and $b_2$. Then $h_9$ is strictly decreasing from $[0,1]$ onto $[\overline{\alpha}B(a,b_2)/B(a,b_1),1]$, with $a_0^{(1)}=h_9(0)=1$,
$$h_9'(0)=a_1^{(1)}=-\frac{a(a+1)(b_2-b_1)}{(c_1+1)(c_2+1)}, ~h_9'(1^-)=-\infty,$$
and for $n\in\IN$,
\begin{align}\label{an(1)}
a_n^{(1)}=\frac{(a)_n(b_1)_n}{n!(c_1+1)_n}-\sum_{k=0}^{n-1}\frac{(a)_{n-k}(b_2)_{n-k}a_k^{(1)}}{(n-k)!(c_2+1)_{n-k}}.
\end{align}

(2) Assume that $a, b_1, b_2\in(0,\infty)$, $c_1=a+b_1$ and $c_2=a+b_2$ with $b_1<b_2$, and for $a\in(1,\infty)$, let
$$\bar{b}=\frac{\sqrt{a^2+4a-1}-\sqrt{a^2-1}}{2}\sqrt{\frac{a+1}{a-1}} ~~(>1).$$
If $a\in(0,1]$, or if $a\in(1,\infty)$ and $b_2\in(0,\bar{b}]$, then $a_n^{(1)}<0$ for all $n\in\IN$ so that the functions $(-h_9')$ and $h_{10}(r)\equiv[1-h_9(r)]/r$ are both absolutely monotone on $(0,1)$, with $h_{10}(0^+)=-a_1^{(1)}$ and $h_{10}(1^-)=1-\overline{\alpha}B(a,b_2)/B(a,b_1)$. In particular, if $b_1=1-a$ and $b_2=3/2-a$ for $a\in(0,1)$, then $h_9$ is strictly decreasing and concave from $[0,1]$ onto $[2\alpha\Gamma(3/2-a)/[\sqrt{\pi}\Gamma(1-a)],1]$, and furthermore, $(-h_9')$ and $h_{10}$ are both absolutely monotone on $(0,1)$ with ranges $(-a_1^{(1)},\infty)$ and $(-a_1^{(1)}, 1-4\Gamma(5/2-a)/[3\sqrt{\pi}\Gamma(2-a)])$, respectively.
\end{lemma}

\begin{proof}
(1) Clearly, $h_9(0)=a_0^{(1)}=1$, and
$$h_9(1)=\frac{\Gamma(b_2+1)\Gamma(c_1+1)}{\Gamma(b_1+1)\Gamma(c_2+1)}=\overline{\alpha}\frac{B(a,b_2)}{B(a,b_1)}$$
by \cite[15.1.20]{AS}. In the particular case when $b_1=1-a$ and $b_2=3/2-a$ for $a\in(0,1)$,
$$h_9(1)=\alpha\frac{\Gamma(3/2-a)}{\Gamma(1-a)\Gamma(3/2)}=\frac{4\Gamma(5/2-a)}{3\sqrt{\pi}\Gamma(2-a)}.$$

By \eqref{Def-2F1}, we have
\begin{align}\label{h9}
h_9(r)=\sum_{n=0}^\infty a_n^{(1)}r^n=\left[\sum_{n=0}^\infty\frac{(a)_n(b_1)_n}{n!(c_1+1)_n}r^n\right]
\left[\sum_{n=0}^\infty\frac{(a)_n(b_2)_n}{n!(c_2+1)_n}r^n\right]^{-1},
\end{align}
from which \eqref{an(1)} follows.

For $n\in\IN_0$, set
$$d_n=\left[\frac{(a)_n(b_1)_n}{n!(c_1+1)_n}\right]\left[\frac{(a)_n(b_2)_n}{n!(c_2+1)_n}\right]^{-1}=\frac{(b_1)_n(c_2+1)_n}{(b_2)_n(c_1+1)_n}.$$
Then for $n\in\IN_0$,
$$\frac{d_{n+1}}{d_n}=\frac{(n+b_1)(n+c_2+1)}{(n+b_2)(n+c_1+1)}=1-\frac{(a+1)(b_2-b_1)}{(n+b_2)(n+c_1+1)}<1.$$
Hence $d_n$ is strictly decreasing in $n\in\IN_0$, and $h_9$ is strictly decreasing on $(0,1)$ by \cite[Lemma 2.1]{PV}.

Differentiation gives
\begin{align}\label{h9'}
h_9'(r)=&\frac{1}{F(a,b_2;c_2+1;r)^2}\left[\frac{ab_1}{c_1+1}F(a+1,b_1+1;c_1+2;r)F(a,b_2;c_2+1;r)\right.\nonumber\\
&\left.-\frac{ab_2}{c_2+1}F(a,b_1;c_1+1;r)F(a+1,b_2+1;c_2+2;r)\right]\nonumber\\
=&-ah_9(r)F(a+1,b_2+1;c_2+2;r)\left[\frac{b_2}{(c_2+1)F(a,b_2;c_2+1;r)}\right.\nonumber\\
&\left.-\frac{b_1F(a+1,b_1+1;c_1+2;r)}{(c_1+1)F(a,b_1;c_1+1;r)F(a+1,b_2+1;c_2+2;r)}\right].
\end{align}
By \cite[15.1.20 \& 15.3.11]{AS},
\begin{align*}
&\lim_{r\to1^-}\left[\frac{b_2}{(c_2+1)F(a,b_2;c_2+1;r)}-\frac{b_1F(a+1,b_1+1;c_1+2;r)}{(c_1+1)F(a,b_1;c_1+1;r)F(a+1,b_2+1;c_2+2;r)}\right]\\
&=\frac{b_2}{(c_2+1)F(a,b_2;c_2+1;1)}-\frac{b_1}{(c_1+1)F(a,b_1;c_1+1;1)}\lim_{r\to1}\frac{F(a+1,b_1+1;c_1+2;r)}{F(a+1,b_2+1;c_2+2;r)}\\
&=\frac{b_2\Gamma(a+1)\Gamma(b_2+1)}{\Gamma(c_2+2)}-\frac{b_1\Gamma(a+1)\Gamma(b_1+1)B(a+1,b_2+1)}{\Gamma(c_1+2)B(a+1,b_1+1)}\\
&=\left(b_2-b_1\right)B(a+1,b_2+1)>0,
\end{align*}
and hence by \eqref{h9'}, we obtain the limiting values
$$h_9'(0)=a_1^{(1)}=-\frac{a(a+1)(b_2-b_1)}{(c_1+1)(c_2+1)}, ~h_9'(1^-)=-\infty,$$
and
\begin{align}
&\lim_{r\to1^-}\sqrt{1-r}h_9'(r)\nonumber\\
=&-ah_9(1)\lim_{r\to1^-}\left\{\sqrt{1-r}F(a+1,b_2+1;c_2+2;r)\left[\frac{b_2}{(c_2+1)F(a,b_2;c_2+1;r)}\right.\right.\nonumber\\
&\left.\left.-\frac{b_1F(a+1,b_1+1;c_1+2;r)}{(c_1+1)F(a,b_1;c_1+1;r)F(a+1,b_2+1;c_2+2;r)}\right]\right\}\nonumber\\
=&-a\left(b_2-b_1\right)B(a+1,b_2+1)h_9(1)\lim_{r\to1^-}\sqrt{1-r}F(a+1,b_2+1;c_2+2;r)=0.\label{h9'(1)}
\end{align}

(2) First, observe that $\bar{b}>1$, which can be easily verified. It follows from \eqref{h9} that
\begin{align}\label{ak(1)}
\frac{(a)_n(b_1)_n}{n!(c_1+1)_n}=\sum_{k=0}^n\frac{(a)_{n-k}(b_2)_{n-k}a_k^{(1)}}{(n-k)!(c_2+1)_{n-k}}.
\end{align}
For $n\in\IN$ and $k=1, 2, ..., n$, let $m=n-k$,
\begin{align*}
S_n^{(1)}=&\frac{n!(c_2+1)_n}{(a)_n(b_2)_n}\cdot\frac{(a)_n(b_1)_n}{n!(c_1+1)_n}=\frac{(b_1)_n(c_2+1)_n}{(b_2)_n(c_1+1)_n},\\
S_n^{(2)}=&\frac{n!(c_2+1)_n}{(a)_n(b_2)_n}\cdot\sum_{k=0}^n\frac{(a)_{n-k}(b_2)_{n-k}a_k^{(1)}}{(n-k)!(c_2+1)_{n-k}}\\
=&1+\frac{n!(c_2+1)_n}{(a)_n(b_2)_n}\sum_{k=1}^n\frac{(a)_{n-k}(b_2)_{n-k}a_k^{(1)}}{(n-k)!(c_2+1)_{n-k}}
\end{align*}
and
\begin{align*}
\Delta_{m,\,n}^{(1)}
=&\left(c_2+1-ab_2\right)(m+n)+2mn+c_2(c_2+1)-ab_2(c_2+2)\\
=&(a+1)(m+n)+2mn+b_2(1-a)\left(m+n+a+b_2+1\right)+a(a+1).
\end{align*}
Then $S_n^{(1)}=d_n=S_n^{(2)}$ by \eqref{an(1)}, so that $S_n^{(2)}$ is strictly decreasing in $n\in\IN$ by the monotonicity of $d_n$. Hence
\begin{align*}
0>&\frac{(a)_{n+1}(b_2)_{n+1}}{n!(c_2+1)_n}\left(S_{n+1}^{(2)}-S_n^{(2)}\right)\\
=&(n+1)(n+c_2+1)\sum_{k=1}^{n+1}\frac{(a)_{m+1}(b_2)_{m+1}a_k^{(1)}}{(m+1)!(c_2+1)_{m+1}}
-(n+a)(n+b_2)\sum_{k=1}^n\frac{(a)_m(b_2)_ma_k^{(1)}}{m!(c_2+1)_m}\\
=&(n+1)(n+c_2+1)a_{n+1}^{(1)}\\
&-\sum_{k=1}^n\frac{(a)_m(b_2)_ma_k^{(1)}}{(m+1)!(c_2+1)_{m+1}}
\left\{\left(n^2+c_2n+ab_2\right)\left[m^2+(c_2+2)m+c_2+1\right]\right.\\
&\left.-\left[n^2+(c_2+2)n+c_2+1\right]\left(m^2+c_2m+ab_2\right)\right\}\\
=&(n+1)(n+c_2+1)a_{n+1}^{(1)}\\
&-\sum_{k=1}^n\frac{(n-m)(a)_m(b_2)_ma_k^{(1)}}{(m+1)!(c_2+1)_{m+1}}
\left\{\left(c_2+1-ab_2\right)(n+m)+2mn+\left[c_2(c_2+1)-ab_2(c_2+2)\right]\right\}\\
=&(n+1)(n+c_2+1)a_{n+1}^{(1)}-\sum_{k=1}^n\frac{k(a)_m(b_2)_ma_k^{(1)}}{(m+1)!(c_2+1)_{m+1}}\Delta_{m,\,n}^{(1)},
\end{align*}
which yields
\begin{align}\label{Ineqan(1)}
a_{n+1}^{(1)}<\frac{1}{(n+1)(n+c_2+1)}\sum_{k=1}^n\frac{k(a)_m(b_2)_ma_k^{(1)}}{(m+1)!(c_2+1)_{m+1}}\Delta_{m,\,n}^{(1)}.
\end{align}

Clearly, $\Delta_{m,\,n}^{(1)}>0$ for $a\in(0,1]$, and if $ab_2\leq c_2(c_2+1)/(c_2+2)$, then
$$\Delta_{m,\,n}^{(1)}\geq\left(c_2+1-ab_2\right)(m+n)>\left(c_2+1-ab_2\right)n>0.$$
On the other hand,
\begin{align*}
ab_2&\leq c_2(c_2+1)/(c_2+2)\Longleftrightarrow ab_2\left(a+b_2+2\right)\leq a^2+2ab_2+b_2^2+a+b_2\\
&\Longleftrightarrow\Delta(a,b_2)\equiv -(a-1)b_2^2-\left(a^2-1\right)b_2+a(a+1)\geq0.
\end{align*}
For each $a\in(1,\infty)$, $\Delta(a,b_2)$ is strictly decreasing in $b_2$ from $(0,\infty)$ onto $(-\infty,a(a+1))$, and hence has the unique zero $\bar{b}$ such that $\Delta(a,b_2)\geq0$ for $a\in(1,\infty)$ and $b_2\in(0,\bar{b}]$. Consequently,
\begin{align}\label{IneqDelmn}
\Delta_{m,\,n}^{(1)}>0 \mbox{~ if ~} a\in(0,1], \mbox{~ or if ~} a\in(1,\infty) \mbox{~ and ~} b_2\in\bigl(0,\bar{b}\bigr].
\end{align}
Since $a_1^{(1)}<0$, it is easy for us to apply \eqref{Ineqan(1)}--\eqref{IneqDelmn} and the mathematical induction to prove that $a_n^{(1)}<0$ for all $n\in\IN$, provided that $a\in(0,1]$, or $a\in(1,\infty)$ and $b_2\in(0,\bar{b}]$. This shows that if $a\in(0,1]$, or if $a\in(1,\infty)$ and $b_2\in(0,\bar{b}]$, then $(-h_9')$ and $h_{10}(r)=-\sum_{n=0}^\infty a_{n+1}^{(1)}r^n$ are both absolutely monotone on $(0,1)$.

The remaining conclusions in part (2) are clear.
\end{proof}

\begin{corollary}\label{Col2.4}
(1) Suppose that $h_{11}(r)=F(1/2,1/2;2;r)/F(1/2,1;5/2;r)=\sum_{k=0}^\infty a_k^{(2)}r^k$ for $r\in(0,1)$.
Then $h_{11}(0^+)=a_0^{(2)}=1$, $h_{11}(1^-)=8/(3\pi)$, and for $k\in\IN$,
\begin{align}\label{Defan}
a_k^{(2)}=\frac{1}{k+1}\left[\frac{(1/2)_k}{k!}\right]^2-3\sum_{j=0}^{k-1}\frac{a_j^{(2)}}{[2(k-j)+1][2(k-j)+3]}<0,
\end{align}
with $a_1^{(2)}=-3/40$, $a_2^{(2)}=-267/11200$, $a_3^{(2)}=-32279/2688000$, $a_4^{(2)}=-40472969/5519360000$ and $a_5^{(2)}=-0.004981513\cdots$. In particular, for each $n\in\IN_0$, the functions $(-h_{11}')$ and
$$h_{12}(r)\equiv\frac{1}{r^{n+1}}\left[\sum_{k=0}^na_k^{(2)}r^k-h_{11}(r)\right]$$
are both absolutely monotone on $(0,1)$ with $h_{12}(0^+)=-a_{n+1}^{(2)}$,
$h_{12}(1)=\tilde{a}_n\equiv\sum_{k=0}^na_k^{(2)}-8/(3\pi)$, $\tilde{a}_4=\sum_{k=0}^4 a_k^{(2)}-8/(3\pi)=0.0353442\cdots$,
and for each $n\in\IN_0$ and all $r\in(0,1)$,
\begin{align}\label{ineqwv'}
\sum_{k=0}^{n+1} a_k^{(2)}r^k-\tilde{a}_{n+1}r^{n+2}<&h_{11}(r)<\sum_{k=0}^{n+2} a_k^{(2)}r^k
\end{align}
and
\begin{align}\label{ineqwv}
\sum_{k=0}^4 a_k^{(2)}r^k-\tilde{a}_4r^5<&h_{11}(r)<\sum_{k=0}^5 a_k^{(2)}r^k.
\end{align}

(2) For $r\in(0,1)$, set $t=\sqrt{r}$, and suppose that
$$h_{13}(r)=\frac{\E(t)-(1-r)\K(t)}{1-(1-r)(\arth\,t)/t}=\sum_{n=0}^\infty a_n^{(3)}r^n.$$
Then $h_{13}(r)\equiv(3\pi/8)h_{11}(r)$, $a_n^{(3)}=(3\pi/8) a_n^{(2)}$ for $n\in\IN_0$, $h_{13}(0^+)=a_0^{(3)}=3\pi/8$, $h_{13}(1^-)=1$, and $a_n^{(3)}<0$ for all $n\in\IN$. In particular, for all $r\in(0,1)$,
\begin{align}\label{ineqh14}
\frac{3\pi}{8}\left(\sum_{k=0}^4 a_k^{(2)}r^{2k}-\tilde{a}_4r^{10}\right)
<\frac{\E(r)-r^{\,\prime2}\K(r)}{1-\bigl(r^{\,\prime2}\arth\,r\bigr)\big{/}r}
<\frac{3\pi}{8}\sum_{k=0}^5 a_k^{(2)}r^{2k}.
\end{align}
\end{corollary}

\begin{proof}
(1) Let $h_9$ and $a_n^{(1)}$ be as in Lemma \ref{lm2.3}. Clearly, $h_{11}(r)=h_9(r)$, $a_n^{(2)}=a_n^{(1)}$ and $h_{12}(r)=-\sum_{k=0}^\infty a_{k+n+1}^{(2)}r^k$ for $a=b_1=b_2/2=1/2$, so that \eqref{Defan} holds. Hence it follows from Lemma \ref{lm2.3} that for $r\in(0,1)$,
$$-a_{n+1}^{(2)}-a_{n+2}^{(2)}r<h_{12}(r)<-a_{n+1}^{(2)}+\left[h_{12}(1^-)+a_{n+1}^{(2)}\right]r.$$
Computation and  \eqref{Defan} give the values of $a_i^{(2)}$ ($1\leq i\leq5$). Consequently, the conclusions in part (1) except for \eqref{ineqwv} follow. Taking $n=4$ in \eqref{ineqwv'}, we obtain \eqref{ineqwv}.

(2) By \eqref{Def-2F1}--\eqref{arth}, it is easy to verify that for $r\in(0,1)$,
\begin{align}\label{Form1}
\frac{\E(r)-r^{\,\prime2}\K(r)}{r^2}=\frac\pi4F\left(\frac12,\frac12;2;r^2\right) \mbox{~and~}
\frac{1-r^{\,\prime2}(\arth\,r)/r}{r^2}=\frac23F\left(\frac12,1;\frac52;r^2\right).
\end{align}
Hence $h_{13}(r)=3\pi h_{11}(r)/8$, and part (2) follows from part (1).
\end{proof}

\section{\normalsize Proof of Theorem \ref{th1.1}}\label{Sec3}

(1) Let $f$ be as in \eqref{IneqA'}. Then by \eqref{Def-2F1} and \eqref{v01ro1}--\eqref{vrn0infty},
\begin{align}
F(a,b;c+1;r)=&\sum_{n=0}^\infty\frac{(a)_n\varphi_n(b)}{n!}r^n, \,F(a,b;c;r)=\sum_{n=0}^\infty\frac{(a)_n\rho_n(b)}{n!}r^n,\label{SerF}\\
f(a,0^+,r)=&\lim_{b\to0^+}\left[\sum_{n=0}^\infty\frac{(a)_n\varphi_n(b)}{n!}r^n\right]
\left[\sum_{n=0}^\infty\frac{(a)_n\rho_n(b)}{n!}r^n\right]^{-1}=1\label{f0}
\end{align}
and
\begin{align}
f(a,\infty,r)=&\lim_{b\to\infty}\left[\sum_{n=0}^\infty\frac{(a)_n\varphi_n(b)}{n!}r^n\right]
\left[\sum_{n=0}^\infty\frac{(a)_n\rho_n(b)}{n!}r^n\right]^{-1}=1.\label{finfty}
\end{align}
By \eqref{f0}--\eqref{finfty}, $f$ is not monotone in $b\in(0,\infty)$, and we can roughly assert that Conjecture \ref{ConjA} is not true.

By \eqref{v01ro1}--\eqref{vrn0infty}, \eqref{tlinf}--\eqref{rl(0)} and Lemma \ref{lm2.1}, and by differentiation, we obtain
\begin{align}\label{f'}
\frac{F(a,b;c;r)^2}{r}\frac{\partial f}{\partial b}
=&\frac1r\left[F(a,b;c;r)\frac{\partial F(a,b;c+1;r)}{\partial b}-F(a,b;c+1;r)\frac{\partial F(a,b;c;r)}{\partial b}\right]\nonumber\\
=&\left[\sum_{n=0}^\infty\frac{(a)_{n+1}\varphi_{n+1}(b)\tau_{n+1}(b)}{(n+1)!}r^n\right]\sum_{n=0}^\infty\frac{(a)_n\rho_n(b)}{n!}r^n\nonumber\\
&-\left[\sum_{n=0}^\infty\frac{(a)_n\varphi_n(b)}{n!}r^n\right]\sum_{n=0}^\infty\frac{(a)_{n+1}\rho_{n+1}(b)\lambda_{n+1}(b)}{(n+1)!}r^n
=\sum_{n=0}^\infty A_nr^n,
\end{align}
where
\begin{align}\label{An}
A_n=&\sum_{k=0}^n\left[\frac{(a)_k\rho_k(b)(a)_{m+1}\varphi_{m+1}(b)\tau_{m+1}(b)}{k!(m+1)!}
-\frac{(a)_{k+1}\rho_{k+1}(b)\lambda_{k+1}(b)(a)_m\varphi_m(b)}{(k+1)!m!}\right]\nonumber\\
=&\sum_{k=0}^n\frac{(a)_k\rho_k(b)(a)_m\varphi_m(b)}{k!m!}
\left[\frac{(m+a)(m+b)}{(m+1)(m+c+1)}\tau_{m+1}(b)-\frac{(k+a)(k+b)}{(k+1)(k+c)}\lambda_{k+1}(b)\right]
\end{align}
and $m=n-k$, with
\begin{align}\label{A0}
A_0=&\frac{a(a+1)}{(c+1)^2}-\frac{a^2}{c^2}=a\frac{b^2-a(a+1)}{c^2(c+1)^2}.
\end{align}

It follows from \eqref{A0} that a necessary condition for the validity of the inequality \eqref{IneqA'} is $A_0\leq0$, namely, $0<b\leq\sqrt{a(a+1)}$.
Consequently, Conjecture \ref{ConjA} is not true if $b>\sqrt{a(a+1)}$.

Next, taking $b_1=1-a$ and $b_2=3/2-a$ for $a\in(0,1)$, we can write \eqref{IneqA1'} as
\begin{align}\label{IneqA2}
f(a,b_1,r)>f(a,b_2,r).
\end{align}
If $a+\sqrt{a(a+1)}<1$, namely, $a\in(0,1/3)$, then $b_2>b_1>\sqrt{a(a+1)}$,
so that by \eqref{f'} and \eqref{A0}, $f$ is strictly increasing in $b\in[b_1,\infty)$ for sufficiently small $r\in(0,1)$. Hence for sufficiently small $r\in(0,1)$,
\begin{align}\label{IneqA2'}
f(a,b_1,r)<f(a,b_2,r),
\end{align}
that is, the inequality in \eqref{IneqA2} is reversed for $a\in(0,1/3)$. This shows that \eqref{IneqA1}, or equivalently, \eqref{IneqA1'} does not hold for $a\in(0,1/3)$ and sufficiently small $r\in(0,1)$.

(2) Clearly, $f_1(a,0)=0$, and
\begin{align*}
&\lim_{r\to1^-}\left[\frac{B(a,1-a)F(a,1-a;2;r)}{B(a,3/2-a)F(a,3/2-a;5/2;r)}-1\right]=\frac{a}{3(1-a)}>0
\end{align*}
by \cite[15.1.20]{AS}. Hence by \cite[15.3.10]{AS},
\begin{align}\label{f1(1)}
f_1(a,1^-)=&\lim_{r\to1^-}\left[\frac{F(a,1-a;2;r)}{B(a,3/2-a)F(a,3/2-a;5/2;r)}\log\frac{e^{R(a,3/2-a)}}{1-r}
-\frac{1}{B(a,1-a)}\log\frac{e^{R(a)}}{1-r}\right]\nonumber\\
=&\frac{F(a,1-a;2;1)}{B(a,3/2-a)F(a,3/2-a;5/2;1)}R(a,3/2-a)-\frac{R(a)}{B(a,1-a)}\nonumber\\
&+\frac{1}{B(a,1-a)}\lim_{r\to1^-}\left[\frac{B(a,1-a)F(a,1-a;2;r)}{B(a,3/2-a)F(a,3/2-a;5/2;r)}-1\right]\log\frac{1}{1-r}=\infty.
\end{align}

Let $h_9$ be as in Lemma \ref{lm2.3} with $b_2=b_1+1/2=3/2-a$ for $a\in(0,1)$. Then by differentiation and \eqref{15.3.3},
\begin{align}\label{f1'}
\frac{\partial f_1}{\partial r}=h_{15}(r)\equiv& h_9(r)\frac{\partial F(a,3/2-a;3/2;r)}{\partial r}+F\left(a,\frac32-a;\frac32;r\right)h_9'(r)-\frac{a(1-a)}{1-r}F(a,1-a;2;r)\nonumber\\
=&\frac{a^2}{3(1-r)}F(a,1-a;2;r)+F\left(a,\frac32-a;\frac32;r\right)h_9'(r).
\end{align}
By \eqref{h9'}, \eqref{f1'} and Lemma \ref{lm2.3}(1), we obtain the limiting values
\begin{align}
h_{15}(0)=&a^2/3-a(a+1)/10=a(7a-3)/30 \label{limh15}
\end{align}
and
\begin{align}
h_{15}(1^-)=&\lim_{r\to1^-}\frac{1}{1-r}\left[\frac{a^2}{3}F(a,1-a;2;r)+\sqrt{1-r}F\left(a,\frac32-a;\frac32;r\right)\cdot\sqrt{1-r}h_9'(r)\right]
=\infty.\label{limh15''}
\end{align}

From \eqref{f1(1)}--\eqref{limh15''}, it is easy for us to obtain the first and second assertions in part (2). By part (1) or by the properties of $f_1$, neither the first (or equivalently, the second) inequality in \eqref{Lam01} nor its reversed one is always valid for all $a,r\in(0,1)$.

(3) Clearly, $f_2(a,0)=\alpha-1=a/[3(1-a)]$. For each $a\in(0,1)$, and for $r\in(0,1)$, let
$$h_{16}(r)=(3-2a)\Gamma\left(\frac32-a\right)F\left(a,\frac32-a;\frac52;r\right)-3\Gamma\left(\frac32\right)\Gamma(2-a)F(a,1-a;2;r).$$
Then by \cite[15.1.20]{AS}, it is easy for us to obtain the limiting value $h_{16}(1^-)=0$,
and hence by l'H\^opital's rule,
\begin{align*}
\lim_{r\to1^-}\frac{h_{16}(r)}{r^{\,\prime}}
=&\lim_{r\to1^-}\left[\frac{3a(1-a)}{2}\Gamma(3/2)\Gamma(2-a)r^{\,\prime}F(a+1,2-a;3;r)\right.\\
&\left.-\frac{a(3-2a)^2}{5}\Gamma(3/2-a)r^{\,\prime}F\left(a+1,\frac52-a;\frac72;r\right)\right]=0.
\end{align*}
Consequently, by \cite[15.1.20 \& 15.3.10]{AS},
\begin{align*}
f_2(a,1^-)=&\lim_{r\to1^-}\left[\frac{\alpha}{B(a,1-a)}\log\frac{e^{R(a)}}{1-r}
-\frac{F(a,1-a;2;r)}{B(a,3/2-a)F(a,3/2-a;5/2;r)}\log\frac{e^{R(a,3/2-a)}}{1-r}\right]\\
=&\frac{\alpha R(a)}{\Gamma(a)\Gamma(1-a)}-\frac{\Gamma(3/2)F(a,1-a;2;1)R(a,3/2-a)}{\Gamma(a)\Gamma(3/2-a)F(a,3/2-a;5/2;1)}\\
&+\lim_{r\to1^-}\left[\frac{\alpha}{\Gamma(a)\Gamma(1-a)}-\frac{\Gamma(3/2)F(a,1-a;2;r)}{\Gamma(a)\Gamma(3/2-a)F(a,3/2-a;5/2;r)}\right]\log\frac{1}{1-r}\\
=&\frac{\alpha R(a)}{\Gamma(a)\Gamma(1-a)}-\frac{\Gamma(3/2)\Gamma(5/2-a)R(a,3/2-a)}{\Gamma(a)\Gamma(3/2-a)\Gamma(2-a)\Gamma(5/2)}
\\
&+\frac{1}{\Gamma(a)}\lim_{r\to1^-}\left[\frac{\alpha}{\Gamma(1-a)}
-\frac{\Gamma(3/2)F(a,1-a;2;r)}{\Gamma(3/2-a)F(a,3/2-a;5/2;r)}\right]\log\frac{1}{1-r}\\
=&\frac{\alpha[R(a)-R(a,3/2-a)]}{\Gamma(a)\Gamma(1-a)}+\frac{\lim\limits_{r\to1^-}\left[\frac{h_{16}(r)}{r^{\,\prime}}\cdot \left(r^{\,\prime}\log\frac{1}{1-r}\right)\right]}{3\Gamma(a)\Gamma(2-a)\Gamma(3/2-a)F(a,3/2-a;5/2;1)}
\\
=&\alpha[\psi(3/2-a)-\psi(1-a)]/B(a,1-a)=\alpha\beta.
\end{align*}

Let $h_9$ be as in \eqref{f1'}. Then by differentiation and \cite[15.3.3]{AS},
\begin{align}\label{g'}
\frac{\partial f_2}{\partial r}
=&a(1-a)\alpha F(a+1,2-a;2;r)-\frac{a(3-2a)}{3}F\left(a+1,\frac52-a;\frac52;r\right)h_9(r)\nonumber\\
&-F\left(a,\frac32-a;\frac32;r\right)h_9'(r)
=F\left(a,\frac32-a;\frac32;r\right)\cdot[-h_9'(r)],
\end{align}
which is a product of two absolutely monotone functions of $r\in(0,1)$ by Lemma \ref{lm2.3}(2). Hence the absolute monotonicity of $f_2$ in $r\in(0,1)$ follows. The monotonicity of $\delta_n$ with respect to $n\in\IN_0$ is clear.

By \eqref{Ka},
$$\frac\pi2\frac{d F\bigr(a,1-a;1;r^2\bigr)}{dr}=\frac{d\K_a(r)}{dr},$$
so that by \eqref{15.3.3} and \cite[Theorem 4.1]{AQVV},
\begin{align*}
\frac{\pi a(1-a)r}{r^{\,\prime 2}} F\left(a,1-a;2;r^2\right)
=2(1-a)\frac{\E_a(r)-r^{\,\prime 2}\K_a(r)}{rr^{\,\prime 2}},
\end{align*}
from which it follows that
\begin{align}\label{dKa}
F(a,1-a;2;r)=\frac{2}{\pi a}\cdot\frac{\E_a(t)-t^{\,\prime2}\K_a(t)}{t^2},
\end{align}
where $t=\sqrt{r}$. Hence
\begin{align}\label{ExprLam1}
\Lambda_1\left(a,r^2\right)=&\frac{2}{\pi a}\cdot\frac{\E_a(r)-r^{\,\prime2}\K_a(r)}{r^2}\cdot\frac{F\bigl(a,3/2-a;3/2;r^2\bigr)}{ F\bigl(a,3/2-a;5/2;r^2\bigr)}
\end{align}
and
\begin{align}\label{Exprf2}
\frac\pi2f_2\left(a,r^2\right)=\alpha\K_a(r)-\frac{\E_a(r)-r^{\,\prime2}\K_a(r)}{ ar^2}\frac{F\bigl(a,3/2-a;3/2;r^2\bigr)}{F\bigl(a,3/2-a;5/2;r^2\bigr)}=\frac\pi2\sum_{k=0}^\infty a_kr^{2k}.
\end{align}
This yields \eqref{Expr1Ka} and \eqref{Expr1K}.

The double inequality \eqref{Ineqf2} follows from \eqref{Exprf2} and the absolute monotonicity of the function
$$r\mapsto \frac{1}{r^{n+1}}\left[f_2(a,r)-\sum_{k=0}^na_kr^k\right]=\sum_{k=0}^\infty a_{k+n+1}r^k.$$
By Lemma \ref{lm2.2}(1), $\beta<1$. The double inequality \eqref{IneqLam1}, and the fourth and fifth inequalities in \eqref{IneqA1''} follow from the monotonicity of $f_2$ in $r\in(0,1)$, and \eqref{IneqLam1'} follows from \eqref{IneqLam1}. The first, second, third and eighth inequalities in \eqref{IneqA1''} are clear.  Since
$$P_{n+1}\left(a,\sqrt{r}\right)=r^{n+2}\sum_{k=0}^{n+2}a_k+\sum_{k=0}^{n+2}a_k\left(r^k-r^{n+2}\right)>r^{n+2}\sum_{k=0}^{n+2}a_k+a_0\left(1-r^{n+2}\right),$$
the sixth inequality in \eqref{IneqA1''} follows. Clearly,
\begin{align*}
r^{n+2}\sum_{k=0}^{n+2}a_k+a_0\bigl(1-r^{n+2}\bigr)>\eta\Longleftrightarrow r^{n+2}\sum_{k=1}^{n+2}a_k>-\frac{a^2}{3(1-a)(3-2a)},
\end{align*}
which is obviously true, and hence the seventh inequality in \eqref{IneqA1''} holds. It is clear that as constant depending only on $a$, each $\alpha$ in \eqref{Ineqf2}--\eqref{IneqLam1} and \eqref{IneqA1''} is best possible.

Next, by \eqref{Def-2F1}, we have the following Maclaurin series
\begin{align}
F(a,1-a;1;r)=&1+\frac{a(1-a)}{\Gamma(a+1)\Gamma(2-a)}\sum_{n=1}^\infty\frac{\Gamma(n+a)\Gamma(n+1-a)}{(n!)^2}r^n,\label{F'a1}\\
F(a,1-a;2;r)=&1+\frac{a(1-a)}{\Gamma(a+1)\Gamma(2-a)}\sum_{n=1}^\infty\frac{\Gamma(n+a)\Gamma(n+1-a)}{(n+1)!n!}r^n,\label{F'a2}\\
F\left(a,\frac32-a;\frac32;r\right)=&1+\frac{a}{\Gamma(a+1)\Gamma(3/2-a)}\sum_{n=1}^\infty\frac{\Gamma(n+a)\Gamma(n+3/2-a)}{n!(3/2)_n}r^n\label{F'a3}
\end{align}
and
\begin{align}
F\left(a,\frac32-a;\frac52;r\right)=&1+\frac{a}{\Gamma(a+1)\Gamma(3/2-a)}\sum_{n=1}^\infty\frac{\Gamma(n+a)\Gamma(n+3/2-a)}{n!(5/2)_n}r^n.\label{F'a4}
\end{align}
It follows from \eqref{F'a1}--\eqref{F'a4} and Lemma \ref{lm2.2}(1) that
$$\Lambda(0^+,r)=1=\alpha(0^+)=\alpha(0^+)[1-\beta(0^+)].$$
Hence each inequality in \eqref{Ineqf2} and \eqref{IneqA1''} is sharp as $a\to0$.

The remaining conclusion in part (3) is clear. ~$\Box$

\section{\normalsize A Counterexample for \eqref{IneqA1} and a Corollary of Theorem \ref{th1.1} }\label{Sec4}

Obviously, by Theorem \ref{th1.1}, one can construct infinitely many counterexamples for \eqref{IneqA} and \eqref{IneqA1}. In this section, we give a detailed counterexample for \eqref{IneqA1}, and derive a corollary of Theorem \ref{th1.1}.

\bigskip
{\it 4.1 ~A Detailed Counterexample for \eqref{IneqA1}}

\bigskip
Take $a=1/4$ in \eqref{IneqA1}, or equivalently, \eqref{IneqA1'}, and let
\begin{align*}
f_5(r)=F\left(\frac14,\frac34;2;r\right)F\left(\frac14,\frac54;\frac32;r\right)-F\left(\frac14,\frac34;1;r\right)
F\left(\frac14,\frac54;\frac52;r\right)
\end{align*}
for $r\in(0,1)$. Then $f_5(0)=0$, and \eqref{IneqA1'} becomes
\begin{align}\label{IneqA4}
\frac{F(1/4,3/4;2;r)}{F(1/4,3/4;1;r)}>\frac{F(1/4,5/4;5/2;r)}{F(1/4,5/4;3/2;r)},
\end{align}
which is valid for all $r\in(0,1)$ if and only if $f_5(r)>0$ for all $r\in(0,1)$.
By \cite[15.3.10]{AS}, it is easy to obtain the limiting value $f_5(1^-)=\infty$, and hence there exists a number $r_6\in(0,1)$ such that
\begin{align}\label{Signf5}
f_5(r)>0 \mbox{~ for ~} r\in[r_6,1).
\end{align}
Differentiation gives
\begin{align}\label{f5'}
f_5^{\,\prime}(r)=&\frac{3}{32}F\left(\frac54,\frac74;3;r\right)F\left(\frac14,\frac54;\frac32;r\right)
+\frac{5}{24}F\left(\frac14,\frac34;2;r\right)F\left(\frac54,\frac94;\frac52;r\right)\nonumber\\
&-\frac{3}{16}F\left(\frac54,\frac74;2;r\right)F\left(\frac14,\frac54;\frac52;r\right)
-\frac18F\left(\frac14,\frac34;1;r\right)F\left(\frac54,\frac94;\frac72;r\right)
\end{align}
with $f_5^{\,\prime}(0)=-1/96$, so that there exists a number $r_7\in(0,1)$ such that $f_5$ is strictly decreasing on $(0,r_7]$ and $f_5(r)<f_5(0)=0$ for $r\in(0,r_7]$. This, together with \eqref{Signf5}, shows that neither the inequality in \eqref{IneqA1'} nor its reversed one always holds for all $r\in(0,1)$, provided that $a=1/4$.

\bigskip
{\it 4.2 ~A Corollary of Theorem \ref{th1.1}}

\bigskip
From Theorem \ref{th1.1}, we can derive the sharp bounds for $\K_a(r)$ and $\K(r)$. Here we give an example as follows. Let $f_2$ and $a_k$ be as in Theorem \ref{th1.1}(3), $h_{13}$ as in Corollary \ref{Col2.4}(2), and suppose that
\begin{align}\label{Deff6}
f_6(r)=\K_a\bigl(\sqrt{r}\bigr)-\frac{\E_a\bigl(\sqrt{r}\bigr)-(1-r)\K_a\bigl(\sqrt{r}\bigr)}{a\alpha r F(a,3/2-a;5/2;r)}F\left(a,\frac32-a;\frac32;r\right)=\sum_{k=0}^\infty a_k^{(4)}r^k
\end{align}
for $a, r\in(0,1)$, where $a_k^{(4)}=a_k^{(4)}(a)$ depends on $a$, and set
\begin{align}\label{f7}
f_7(r)=f_6(r)|_{a=1/2}=\K\left(\sqrt{r}\right)-h_{13}(r)\frac{{\rm{arth}}\sqrt{r}}{\sqrt{r}}.
\end{align}
Then $f_6(r)=\pi f_2(a,r)/(2\alpha)$ by \eqref{Exprf2}, so that by Theorem \ref{th1.1}(3), $f_6$ and $f_7$ are both absolutely monotone on $(0,1)$ with $a_k^{(4)}(a)=\pi a_k/(2\alpha)$, $f_6(0)=\pi\eta/2$, $f_6(1^-)=\pi\beta/2$, $f_7(0^+)=\pi/8$ and $f_7(1^-)=\log2$.
Hence by Theorem \ref{th1.1}(3), \eqref{arth} and Corollary \ref{Col2.4}, we obtain the following corollary immediately.
\begin{corollary}\label{Col4.1}
Let $h_{13}$, $f_6$ and $f_7$ be as Corollary \ref{Col2.4}(2), \eqref{Deff6} and \eqref{f7}, respectively. Then for each $a\in(0,1)$, $f_6$ and $f_7$ are both absolutely monotone on $(0,1)$ with ranges $(\pi\eta/2, \pi\beta/2)$ and $(\pi/8,\log2)$, respectively. In particular, for all $r\in(0,1)$,
\begin{align}\label{ineq2g2}
&\frac\pi8+\frac{3\pi}{8}\left[1-\frac{3}{40}r^2-\frac{267}{11200}r^4-\frac{32279}{2688000}r^6
-\left(\frac{2390041}{2688000}-\frac{8}{3\pi}\right)r^8\right]\frac{{\rm{arth}}\,r}{r}\nonumber\\
&<\frac\pi8+h_{13}\left(r^2\right)\frac{{\rm{arth}}\,r}{r}<\K(r)<\frac\pi8+h_{13}\left(r^2\right)\frac{{\rm{arth}}\,r}{r}+\left(\log2-\frac\pi8\right)r^2\nonumber\\
&<\frac\pi8+\left(\log2-\frac\pi8\right)r^2+\frac{3\pi}{8}\left(1-\frac{3}{40}r^2-\frac{267}{11200}r^4
-\frac{32279}{2688000}r^6-\frac{40472969}{5519360000}r^8\right)\frac{{\rm{arth}}\,r}{r},
\end{align}
in which all the constants are best possible.
\end{corollary}

\section{\normalsize Proof of Theorem \ref{th1.2}}\label{Sec5}

(1) Suppose that $\Lambda_3(a,b_1,b_2,r)=\sum_{n=0}^\infty B_nr^n$. Then
\begin{align*}
\left[\sum_{n=0}^\infty\frac{(a)_n(b_1)_n}{n!(c_1+1)_n}r^n\right]\sum_{n=0}^\infty\frac{(a)_n(b_2)_n}{n!(c_2)_n}r^n=\left(\sum_{n=0}^\infty B_nr^n\right)\sum_{n=0}^\infty\frac{(a)_n(b_2)_n}{n!(c_2+1)_n}r^n,
\end{align*}
by which we obtain $B_0=1$ and
\begin{align}\label{Bn}
B_n=\frac{(a)_n(b_1)_n}{n!(c_1+1)_n}+\sum_{k=0}^{n-1}\frac{(a)_{n-k}(b_2)_{n-k}}{(n-k)!(c_2)_{n-k+1}}
\left[\left(n-k+c_2\right)\frac{(a)_k(b_1)_k}{k!(c_1+1)_k}-c_2B_k\right]
\end{align}
for $n\in\IN$, with $B_1=a\bigl\{b_1/(c_1+1)+b_2/[c_2(c_2+1)]\bigl\}$. Hence
\begin{align}\label{Exprf3}
\frac{f_3(a,b_1,b_2,r)}{r}=&\frac{\Lambda_3(a,b_1,b_2,r)-F(a,b_1;c_1;r)}{r}
=\sum_{n=1}^\infty\left[B_n-\frac{(a)_n(b_1)_n}{n!(c_1)_n}\right]r^{n-1}\nonumber\\
=&\frac{a(b_2-b_1)\bigl(a^2+a-b_1b_2\bigr)}{c_1c_2(c_1+1)(c_2+1)}
+\sum_{n=1}^\infty\left[B_{n+1}-\frac{(a)_{n+1}(b_1)_{n+1}}{(n+1)!(c_1)_{n+1}}\right]r^n
\end{align}
by which we obtain the limiting value
\begin{align*}
\lim_{r\to0}\frac{\Lambda_3(a,b_1,b_2,r)-F(a,b_1;c_1;r)}{r}=\frac{a(b_2-b_1)\bigl[a(a+1)-b_1b_2\bigr]}{c_1c_2(c_1+1)(c_2+1)}.
\end{align*}
This shows that if $b_1b_2>a(a+1)$, then there exists a number $r_0\in(0,1)$ such that
$$\Lambda_3(a,b_1,b_2,r)<F(a,b_1;c_1;r)$$
for $r\in(0,r_0]$. Consequently, we obtain part (1) by \eqref{Lam23}.

(2) Clearly, $f_3(a,b_1,b_2,0)=0$. By \cite[15.1.20 \& 15.3.10]{AS}, we obtain
\begin{align}\label{f3(1)}
f_3(a,b_1,b_2,1^-)=&\lim_{r\to1^-}\left[\frac{F(a,b_1;c_1+1;r)}{B(a,b_2)F(a,b_2;c_2+1;r)}\log\frac{e^{R(a,b_2)}}{1-r}
-\frac{1}{B(a,b_1)}\log\frac{e^{R(a,b_1)}}{1-r}\right]\nonumber\\
=&\frac{R(a,b_2)F(a,b_1;c_1+1;1)}{B(a,b_2)F(a,b_2;c_2+1;r)}-\frac{R(a,b_1)}{B(a,b_1)}\nonumber\\
&+\frac{1}{B(a,b_1)}\lim_{r\to1^-}\left[\frac{B(a,b_1)F(a,b_1;c_1+1;r)}{B(a,b_2)F(a,b_2;c_2+1;r)}-1\right]\log\frac{1}{1-r}=\infty
\end{align}
since by \cite[15.1.20]{AS},
\begin{align*}
\frac{B(a,b_1)F(a,b_1;c_1+1;1)}{B(a,b_2)F(a,b_2;c_2+1;1)}-1=\overline{\alpha}-1=a\frac{b_2-b_1}{b_1c_2}>0.
\end{align*}

Let $h_9$ be as in Lemma \ref{lm2.3}(1). Then $f_3(a,b_1,b_2,r)=F(a,b_2;c_2;r)h_9(r)-F(a,b_1;c_1;r)$, and
\begin{align}\label{f3'}
\frac{\partial f_3}{\partial r}=&f_8(r)\equiv\frac{a^2(b_2-b_1)}{c_1c_2(1-r)}F(a,b_1;c_1+1;r)+F(a,b_2;c_2;r)h_9'(r).
\end{align}
By \eqref{f3'}, Lemma \ref{lm2.3}(1) and \eqref{h9'(1)}, we obtain the limiting values
\begin{align}\label{f8(0)}
f_8(0)=&\frac{a^2(b_2-b_1)}{c_1c_2}+h_9'(0)
=\frac{a(b_2-b_1)[a(a+1)-b_1b_2]}{c_1c_2(c_1+1)(c_2+1)}
\end{align}
and
\begin{align}\label{f8(1)}
f_8(1^-)=&\lim_{r\to1^-}\frac{1}{1-r}\left[\frac{a^2(b_2-b_1)}{c_1c_2}F(a,b_1;c_1+1;r)
+\sqrt{1-r}F(a,b_2;c_2;r)\cdot\sqrt{1-r}h_9'(r)\right]=\infty.
\end{align}

Consequently, by  \eqref{f3(1)}--\eqref{f8(1)}, if $b_1b_2>a(a+1)$ ($b_1b_2<a(a+1)$), then there exists a number $r_4\in(0,1)$ ($r_5\in(0,1)$) such that $f_3$ is strictly decreasing (increasing) in $r\in(0,r_4)$ ($r\in(0,r_5)$, respectively), and there exists a number $r_6\in(0,1)$ such that $f_3$ is positive and strictly increasing in $r\in(r_6,1)$. This yields the first and second assertions in part (2).

For arbitrarily given $a,b_1,b_2\in(0,\infty)$ with $b_1<b_2$, and for $r\in(0,1)$, let $g_1(r)=\Lambda_2(a,b_1,b_2,r)$.
Then $g_1(0)=1$ and $g_1(1^-)=\overline{\alpha}$ by \eqref{limLam}. By differentiation,
\begin{align}\label{g1'}
&(1-r)\left[F(a,b_1;c_1;r)F(a,b_2;c_2+1;r)\right]^2g_1'(r)\nonumber\\
=&F(a,b_1;c_1;r)F(a,b_2;c_2+1;r)\left[\frac{ab_2}{c_2}F(a,b_2;c_2+1;r)F(a,b_1;c_1+1;r)\right.\nonumber\\
&\left.+\frac{ab_1(1-r)}{c_1+1}F(a,b_2;c_2;r)F(a+1,b_1+1;c_1+2;r)\right]\nonumber\\
&-F(a,b_2;c_2;r)F(a,b_1;c_1+1;r)\left[\frac{ab_1}{c_1}F(a,b_1;c_1+1;r)F(a,b_2;c_2+1;r)\right.\nonumber\\
&\left.+\frac{ab_2(1-r)}{c_2+1}F(a,b_1;c_1;r)F(a+1,b_2+1;c_2+2;r)\right],
\end{align}
by which we obtain the limiting value
\begin{align*}
g_1'(0)=\frac{ab_2}{c_2}+\frac{ab_1}{c_1+1}-\frac{ab_1}{c_1}-\frac{ab_2}{c_2+1}
=\frac{a(b_2-b_1)[a(a+1)-b_1b_2]}{c_1c_2(c_1+1)(c_2+1)}.
\end{align*}
Hence if $b_1b_2<a(a+1)$ ($b_1b_2>a(a+1)$), then there exists a number $r_8\in(0,1)$ ($r_9\in(0,1)$) such that $g_1$ is strictly increasing (decreasing) on $(0,r_8]$ ($(0,r_9]$, respectively). Consequently,
$$g_1(r)>g_1(0)=1 ~~(g_1(r)<g_1(0)=1)$$
for $r\in(0,r_8]$ ($r\in(0,r_9]$, respectively). This yields the last conclusion in part (2).

(3) Clearly, $f_4(a,b_1,b_2,r)|_{b_1=b_2}\equiv0$ and $f_4(a,b_1,b_2,0)=\overline{\alpha}-1=a(b_2-b_1)/(b_1c_2)$. Put
\begin{align*}
g_2(r)=&b_1c_2B(a,b_1)F(a,b_2;c_2+1;r)\left[h_9(r)-\overline{\alpha}\frac{B(a,b_2)}{B(a,b_1)}\right].
\end{align*}
Then $g_2(1^-)=0$ by \cite[15.1.20]{AS} and Lemma \ref{lm2.3}(1). By l'H\^opital's rule, and by \cite[15.1.20]{AS} and \eqref{h9'(1)},
\begin{align*}
\lim_{r\to1^-}\frac{g_2(r)}{r^{\,\prime}}=&-b_1c_2B(a,b_1)\lim_{r\to 1^-}\left[\frac{r^{\,\prime}}{r}F(a,b_2;c_2+1;1)h_9'(r)\right]=0.
\end{align*}
Hence by \cite[15.1.20 \& 15.3.10]{AS}, and we obtain the limiting value
\begin{align*}
f_4(a,b_1,b_2,1^-)=&\lim_{r\to1^-}\left[\frac{\overline{\alpha}}{B(a,b_1)}\log\frac{e^{R(a,b_1)}}{1-r}
-\frac{F(a,b_1;c_1+1;r)}{B(a,b_2)F(a,b_2;c_2+1;r)}\log\frac{e^{R(a,b_2)}}{1-r}\right]\\
=&\frac{\overline{\alpha}}{B(a,b_1)}R(a,b_1)-\frac{F(a,b_1;c_1+1;1)}{B(a,b_2)F(a,b_2;c_2+1;1)}R(a,b_2)\\
&-\lim_{r\to1^-}\frac{1}{b_1c_2B(a,b_1)B(a,b_2)F(a,b_2;c_2+1;1)}\lim_{r\to1^-}\left[\frac{g_2(r)}{r^{\,\prime}}\cdot r^{\,\prime}\log\frac{1}{1-r}\right]\\
=&\frac{\overline{\alpha}\Gamma(c_1)}{\Gamma(a)\Gamma(b_1)}\left[R(a,b_1)-R(a,b_2)\right]
=\overline{\alpha}\frac{\psi(b_2)-\psi(b_1)}{B(a,b_1)}=\overline{\alpha}\overline{\beta}.
\end{align*}

Differentiation and Lemma \ref{lm2.3}(1) give
\begin{align}\label{g2'}
\frac{\partial f_4}{\partial r}=&F(a,b_2;c_2;r)[-h_9'(r)]>0
\end{align}
for all $r\in(0,1)$, which yields the monotonicity of $f_4$ in $r\in(0,1)$.

The remaining conclusions in part (3) are clear.

(4) Assume that $a\in(0,1]$, or $a\in(1,\infty)$ and $b_2\in(0,\bar{b}]$. Then it follows from Lemma \ref{lm2.3}(2) and \eqref{g2'} that as a function of $r\in(0,1)$,
$$\frac{\partial f_4}{\partial r}=F(a,b_2;c_2;r)[-h_9'(r)]$$
is a product of two absolutely monotone functions on $(0,1)$. Hence $f_4$ is absolutely monotone in $r\in(0,1)$, and so is the function
$$h_{17}(r)\equiv\frac{1}{r^{n+1}}\left[f_4(a,b_1,b_2,r)-\sum_{k=0}^n \tilde{b}_kr^k\right]=\sum_{k=0}^\infty \tilde{b}_{k+n+1}r^k$$
for each $n\in\IN_0$, with $h_{17}(0^+)=\tilde{b}_{n+1}$ and $h_{17}(1^-)=\bar{\alpha}\bar{\beta}-\sum_{k=0}^n \tilde{b}_k$. By the property of $h_{17}$, we obtain
\begin{align*}
\tilde{b}_{n+1}+\tilde{b}_{n+2}r<h_{17}(r)<h_{17}(0^+)+\left[h_{17}(1^-)-h_{17}(0^+)\right]r=\tilde{b}_{n+1}+\eta_n r,
\end{align*}
which yields \eqref{IneqLam23}. The second and third inequalities in \eqref{IneqLam23'} follow from \eqref{IneqLam23}, and the first and fourth  inequalities in \eqref{IneqLam23'} are clear.  ~$\Box$

\section{\normalsize Concluding Remarks}\label{Sec6}

In this section, we explain that we can disprove Conjecture \ref{ConjA} with other methods, and put forward three open problems and a conjecture.

1. Suppose that $\Lambda_1(a,r)=\sum_{n=0}^\infty C_nr^n$. Then $C_0=1$ and
$$\left[\sum_{n=0}^\infty\frac{(a)_n(1-a)_n}{(n+1)!n!}r^n\right]\sum_{n=0}^\infty\frac{(a)_n(3/2-a)_n}{n!(3/2)_n}r^n=\left(\sum_{n=0}^\infty C_nr^n\right)\sum_{n=0}^\infty\frac{(a)_n(3/2-a)_n}{n!(5/2)_n}r^n,$$
from which we obtain
\begin{align}\label{Cn}
C_n=\frac{(a)_n(1-a)_n}{(n+1)!n!}+\sum_{k=0}^{n-1}\frac{(a)_{n-k}(3/2-a)_{n-k}}{(n-k)!(3/2)_{n-k+1}}\left[\left(n-k+\frac32\right)\frac{(a)_k(1-a)_k}{(k+1)!k!}-\frac32C_k\right]
\end{align}
for $n\in\IN$, with $C_1=a(27-23a)/30$. Hence
\begin{align}\label{Serf1}
\frac{f_1(a,r)}{r}=&\sum_{n=1}^\infty\left[C_n-\frac{(a)_n(1-a)_n}{(n!)^2}\right]r^{n-1}
=\frac{a(7a-3)}{30}+\sum_{n=2}^\infty\left[C_n-\frac{(a)_n(1-a)_n}{(n!)^2}\right]r^{n-1}
\end{align}
and by Theorem \ref{th1.1},
\begin{align}\label{Serf2}
f_2(a,r)=&\sum_{n=0}^\infty a_nr^n=\sum_{n=0}^\infty\left[\alpha\frac{(a)_n(1-a)_n}{(n!)^2}-C_n\right]r^n\nonumber\\
=&\frac{a}{3(1-a)}+\frac{a(1+a)}{10}r+\sum_{n=2}^\infty\left[\alpha\frac{(a)_n(1-a)_n}{(n!)^2}-C_n\right]r^n.
\end{align}
By \eqref{Serf1}, we see that \eqref{Lam01} is not always valid.

2. It follows from \cite[15.3.1 \& 15.3.10]{AS}, \eqref{F'a1}--\eqref{F'a4} and \eqref{arth} that for each $r\in(0,1)$,
\begin{align}
&\Lambda(0^+,r)=1, \mbox{~ and ~} \Lambda(1^-,r)=\frac{F(1/2,1;3/2;r)}{F(1/2,1;5/2;r)}
=\frac{2\bigl(\arth\sqrt{r}\bigr)\big{/}\sqrt{r}}{3\bigl[1-(1-r)\bigl(\arth\sqrt{r}\bigr)\big{/}\sqrt{r}
\bigr]}\label{lim2Lam}
\end{align}
which is strictly increasing from $(0,1)$ onto $(1,\infty)$.

Suppose that $\Lambda(a,r)=\sum_{n=0}^\infty D_nr^n$. Then $D_0=1$ and
$$\Lambda_1(a,r)=F(a,1-a;1;r)\sum_{n=0}^\infty D_nr^n$$
from which it follows that
\begin{align}\label{Dn}
D_n=C_n-\sum_{k=0}^{n-1}D_k\frac{(a)_{n-k}(1-a)_{n-k}}{[(n-k)!]^2}
\end{align}
for $n\in\IN$, with $D_1=a(7a-3)/30$. Hence
\begin{align*}
\frac{\Lambda(a,r)-1}{r}=\frac{a(7a-3)}{30}+\sum_{n=1}^\infty\left\{C_{n+1}-\sum_{k=0}^n\frac{(a)_{n-k+1}(1-a)_{n-k+1}D_k}{[(n-k+1)!]^2}\right\}r^n,
\end{align*}
by which we see that \eqref{Lam01} is not always valid for $r\in(0,1)$ if $a\in(0,3/7)$.

On the other hand,
$$\frac{\partial \Lambda(a,r)}{\partial r}=\sum_{n=0}^\infty (n+1)D_{n+1}r^n,$$
and
$$\lim_{r\to0}\frac{\partial \Lambda(a,r)}{\partial r}=D_1=\frac{a(7a-3)}{30}.$$
Hence by \eqref{limLam}, if $a\in(0,3/7)$, then $\Lambda(a,r)<1$ for sufficiently small $r$. This also shows that \eqref{IneqA1'} is not valid for $r$ close to $0$, provided that $a\in(0,3/7)$.

3. By the series \eqref{F'a1}--\eqref{F'a4}, we obtain the following limiting values
\begin{align}
\lim_{a\to0^+}F(a,1-a;1;r)=&\lim_{a\to0^+}F(a,1-a;2;r)\nonumber\\
=\lim_{a\to0^+}F(a,3/2-a;3/2;r)=&\lim_{a\to0^+}F(a,3/2-a;5/2;r)=1,\label{F'0}
\end{align}
and by differentiation and \cite[6.1.15 \& 6.3.5]{AS}, we have
\begin{align}
\frac{\partial F(a,1-a;2;r)}{\partial a}
=&\sum_{n=1}^\infty\frac{\Gamma(n+a)\Gamma(n+1-a)\psi_{1,\,n}(a)}{\Gamma(a+1)\Gamma(2-a)(n+1)!n!}r^n,\label{F'2}\\
\frac{\partial F(a,3/2-a;3/2;r)}{\partial a}
=&\sum_{n=1}^\infty\frac{\Gamma(n+a)\Gamma(n+3/2-a)\psi_{2,\,n}(a)}{\Gamma(a+1)\Gamma(3/2-a)(3/2)_nn!}r^n\label{F'3}
\end{align}
and
\begin{align}
\frac{\partial F(a,3/2-a;5/2;r)}{\partial a}
=&\sum_{n=1}^\infty\frac{\Gamma(n+a)\Gamma(n+3/2-a)\psi_{2,\,n}(a)}{\Gamma(a+1)\Gamma(3/2-a)(5/2)_nn!}r^n,\label{F'4}
\end{align}
where
\begin{align*}
\psi_{1,\,n}(a)=&a(1-a)[\psi(n+a)-\psi(n+1-a)-\psi(a+1)+\psi(2-a)]+1-2a,\\
\psi_{2,\,n}(a)=&a[\psi(n+a)-\psi(n+3/2-a)-\psi(a+1)+\psi(3/2-a)]+1.
\end{align*}
Applying \eqref{partialF} and \eqref{F'2}--\eqref{F'4}, we obtain
\begin{align}
\lim_{a\to0^+}\frac{\partial F(a,1-a;1;r)}{\partial a}&=\lim_{a\to0^+}\frac{\partial F(a,3/2-a;3/2;r)}{\partial a}=\sum_{n=1}^\infty\frac{r^n}{n}=\log\frac{1}{1-r},\label{dF01}\\
\lim_{a\to0^+}\frac{\partial F(a,1-a;2;r)}{\partial a}&=\sum_{n=1}^\infty\frac{r^n}{n(n+1)}=1-\frac{1-r}{r}\log\frac{1}{1-r}\label{dF01'}
\end{align}
and
\begin{align}
&\lim_{a\to0^+}\frac{\partial F(a,3/2-a;5/2;r)}{\partial a}=\sum_{n=1}^\infty\frac{3 r^n}{n(2n+3)}
=\frac23+\log\frac{1}{1-r}-\frac{2}{r}\left(\frac{\arth\sqrt{r}}{\sqrt{r}}-1\right).\label{dF02}
\end{align}

Logarithmic differentiation gives
\begin{align}\label{Lam'a}
\frac{1}{\Lambda(a,r)}\frac{\partial\Lambda(a,r)}{\partial a}=&\frac{1}{F(a,1-a;2;r)}\frac{\partial F(a,1-a;2;r)}{\partial a}\nonumber\\
&+\frac{1}{F(a,3/2-a;3/2;r)}\frac{\partial F(a,3/2-a;3/2;r)}{\partial a}\nonumber\\
&-\frac{1}{F(a,1-a;1;r)}\frac{\partial F(a,1-a;1;r)}{\partial a}\nonumber\\
&-\frac{1}{F(a,3/2-a;5/2;r)}\frac{\partial F(a,3/2-a;5/2;r)}{\partial a}.
\end{align}
By \eqref{F'0} and \eqref{dF01}--\eqref{Lam'a}, we obtain
\begin{align*}
\lim_{a\to0^+}\frac{\partial\Lambda(a,r)}{\partial a}=&\lim_{a\to0^+}\left[\frac{\partial F(a,1-a;2;r)}{\partial a}+\frac{\partial F(a,3/2-a;3/2;r)}{\partial a}\right.\\
&\left.-\frac{\partial F(a,1-a;1;r)}{\partial a}
-\frac{\partial F(a,3/2-a;5/2;r)}{\partial a}\right]\\
=&\sum_{n=1}^\infty\frac{r^n}{n(n+1)}-3\sum_{n=1}^\infty\frac{r^n}{n(2n+3)}
=-\sum_{n=1}^\infty\frac{r^n}{(n+1)(2n+3)}\\
=&\frac13+\frac{1}{r}\left[2\left(\frac{\arth\sqrt{r}}{\sqrt{r}}-1\right)-\log\frac{1}{1-r}\right],
\end{align*}
which shows that for each $r\in(0,1)$, there exists a number $\check{a}\in(0,\infty)$ such that $\Lambda(a,r)$ is strictly decreasing in $a\in(0,\check{a}]$, and $\Lambda(a,r)<\Lambda(0^+,r)=1$ for $a\in(0,\check{a}]$ by \eqref{lim2Lam}. This also shows that \eqref{IneqA1'} does not hold for $a\in(0,\check{a}]$, and we see that $\Lambda(a,r)$ is not monotone in $a\in(0,1)$ by \eqref{lim2Lam}.

4. By applying \eqref{partialF}, \eqref{F'2}--\eqref{F'4} and \eqref{Lam'a}, it is not difficult for us to prove that for each $r\in(0,1)$,
$$\lim_{a\to1^-}\frac{\partial\Lambda(a,r)}{\partial a}>0,$$
showing that for each $r\in(0,1)$, there exists $\hat{a}\in(0,1)$ such that $\Lambda(a,r)$ is strictly increasing in $a\in[\hat{a},1)$.

It seems that as a function of $a$, $\Lambda(a,r)$ is strictly decreasing and then increasing on $(0,1)$.

5. For $a, b_1, b_2\in(0,\infty)$ and $r\in(0,1)$ with $b_1<b_2$, suppose that $\Lambda_2(a,b_1,b_2,r)=\sum_{n=0}^\infty\overline{B}_nr^n$, and let $B_n$ be as in the proof of Theorem \ref{th1.2}(1). Then $\overline{B}_0=1$ and
$$\sum_{n=0}^\infty B_nr^n=\left(\sum_{n=0}^\infty\overline{B}_nr^n\right)\sum_{n=0}^\infty\frac{(a)_n(b_1)_n}{(c_1)_nn!}r^n,$$
from which it follows that
$$\overline{B}_n=B_n-\sum_{k=0}^{n-1}\frac{(a)_{n-k}(b_1)_{n-k}}{(c_1)_{n-k}(n-k)!}\overline{B}_k$$
for $n\in\IN$, with $\overline{B}_1=a(b_2-b_1)\bigl(a^2+a-b_1b_2\bigr)\big{/}\bigl[c_1c_2(c_1+1)(c_2+1)\bigr]$. Hence
\begin{align}\label{SerLam2}
\frac{\Lambda_2(a,b_1,b_2,r)-1}{r}=\frac{a(b_2-b_1)[a(a+1)-b_1b_2]}{c_1c_2(c_1+1)(c_2+1)}+\sum_{n=1}^\infty\overline{B}_{n+1}r^n
\end{align}
and
\begin{align*}
\lim_{r\to0}\frac{\Lambda_2(a,b_1,b_2,r)-1}{r}=\frac{a(b_2-b_1)[a(a+1)-b_1b_2]}{c_1c_2(c_1+1)(c_2+1)}.
\end{align*}
Consequently, if $b_1b_2>a(a+1)$, then $\Lambda_2(a,b_1,b_2,r)<1$ for sufficiently small $r\in(0,1)$, namely, \eqref{IneqA} in Conjecture \ref{ConjA} is not always valid by \eqref{Lam23}.

6. For $a, b\in(0,\infty)$ and $n\in\IN_0$, let $f$ be as in \eqref{IneqA'}. Then
\begin{align*}
f(a,b,r)=\left[\sum_{n=0}^\infty\frac{(a)_n\varphi_n(b)}{n!}r^n\right]\left[\sum_{n=0}^\infty\frac{(a)_n\rho_n(b)}{n!}r^n\right]^{-1}.
\end{align*}
Applying \eqref{vrn+1} and \cite[Lemma 2.1]{PV}, we can easily show that $f$ is strictly decreasing in $r\in(0,1)$, so that   $f(a,b,r)<f(a,b,0^+)=1=f(a,\infty,r)$ by \eqref{finfty}. Consequently, there is at least a sufficiently large $\tilde{b}>0$ such that $f$ is not decreasing in $b\in[\tilde{b},\infty)$, that is, \eqref{IneqA} does not hold for $b\in[\tilde{b},\infty)$.

7. Based on our studies, the following problems are natural.
\begin{open problem}\label{OP}
(1) For what values of the positive parameters $a$ and $b$, \eqref{IneqA} is valid for all $r\in(0,1)$? Can we find the conditions satisfied by the positive parameters a and b, such that these conditions are sufficient and necessary for \eqref{IneqA} to be valid for all $r\in(0,1)$?

(2) For what values of the parameter $a\in(0,1)$, \eqref{IneqA1} is valid for all $r\in(0,1)$?

(3) In \eqref{IneqA1''} and \eqref{Ineq2g1}, what are the best possible lower bounds, which are independent of $r$, of $\Lambda(a,r)$ and $\Lambda_2(a,b_1,b_2,r)$?
\end{open problem}

In addition, our computation seems to show that the following conjecture is true.
\begin{conjecture}\label{Qcon}
The inequality \eqref{IneqA} holds for all $r\in(0,1)$ and $a\in(0,\infty)$ if and only if $b\in\bigl(0,\sqrt{a(a+1)}\bigr]$, and
the inequality \eqref{IneqA1} holds for all $r\in(0,1)$ if and only if $a\in[3/7,1)$.
\end{conjecture}

\bigskip
{\bf Acknowledgements.} This work is supported by Zhejiang Sci-Tech University through the Fund of The First-Class Discipline Construction Project of Zhejiang Province, China.

\bigskip
{\bf Competing Interests}: The authors declare that they have no conflict of interest.

\bigskip
\noindent Authors' address:

\smallskip
\noindent Song-Liang Qiu:\\
Department of Mathematics,
Zhejiang Sci-Tech University\\
Hangzhou 310018, Zhejiang, China\\
Email address: sl$\_$qiu@zstu.edu.cn

\bigskip
\noindent Xiao-Yan Ma:\\
Department of Mathematics,
Zhejiang Sci-Tech University\\
Hangzhou 310018, Zhejiang, China\\
E-mail address: mxy@zstu.edu.cn

\bigskip
\noindent Xue-Yan Xiang:\\
Department of Mathematics, Lishui University\\
Lishui 323000, Zhejiang, China\\
E-mail address: xxy$\_$81917@126.com


\end{document}